\documentclass[12pt]{amsart}
\usepackage{amssymb}

\usepackage{amscd}
\usepackage{graphics}
\usepackage{amsmath}
\usepackage{epsfig, psfrag, graphics}

\newtheorem{prop}{Proposition}[section]
\newtheorem{thm}[prop]{Theorem}
\newtheorem{lem}[prop]{Lemma}

\theoremstyle{definition}
\newtheorem{Ack}[prop]{Acknowledgments}

\newtheorem{theorem}[prop]{Theorem}
\newtheorem{proposition}[prop]{Proposition}
\newtheorem{lemma}[prop]{Lemma}
\newtheorem{corollary}[prop]{Corollary}
\newtheorem{remark}[prop]{Remark}

\newtheorem{definition}[prop]{Definition}

\input tcilatex

\begin{document}
\author{Emmanuel Breuillard}
\date{March 2008}
\address{Emmanuel Breuillard, Ecole Polytechnique, France}
\email{emmanuel.breuillard@math.polytechnique.fr}

%
%
\title[A strong Tits alternative]
{A strong Tits alternative}

\begin{abstract}
We show that for every integer $d\in \Bbb{N}$, there is $N(d)\in \Bbb{N}$
 such that if $K$ is any field and $F$ is a finite subset of $GL_d(K)$, which
generates a non amenable subgroup, then $F^{N(d)}$ contains two elements, 
which freely generate a non abelian free subgroup. This improves the original 
statement of the Tits alternative. It also implies a growth gap and a co-growth gap for
 non-amenable linear groups, and has consequences about the girth and uniform expansion of small sets in finite
 subgroups of $GL_d(\Bbb{F}_q)$ as well as other diophantine properties of non-discrete subgroups of
 Lie groups.
\end{abstract}
\maketitle%

\section{Introduction}

The goal of this paper is to show the following theorem and some
consequences of it.

\begin{theorem}
\label{main}For every $d\in \Bbb{N}$ there is $N(d)\in \Bbb{N}$ such that if 
$K$ is any field and $F$ a finite symmetric subset of $GL_{d}(K)$ containing 
$1$, either $F^{N(d)}$ contains two elements which freely generate a non
abelian free group, or the group generated by $F$ is virtually solvable
(i.e. contains a finite index solvable subgroup).
\end{theorem}

By $F^{N(d)}=F\cdot ...\cdot F$ we mean the set of elements which can be
written as a product of at most $N(d)$ elements from $F,$ and by symmetric
we mean that if $f\in F$ then $f^{-1}\in F.$ This statement is a
strengthening of the classical Tits alternative \cite{Tits}, which asserts
that any finitely generated subgroup $\left\langle F\right\rangle $ of $%
GL_{d}(K),$ where $K$ is any field, either contains a non abelian free
subgroup or contains a solvable subgroup of finite index. It also improves
earlier strengthenings of the Tits alternative, due to Eskin-Mozes-Oh \cite
{EMO} (for free semigroups) and to T. Gelander and the author \cite{uti}
(for free groups), which showed a statement of a similar form, except that
the integer $N(d)$ depended on the group $\Gamma $ generated by $F$ (not on
the generating set) but was not independent of the field of coefficients.
Note that $N(d)$ cannot be bounded uniformly in $d$ (see Remark \ref{GrH}).

The present paper essentially contains the geometric part of the proof of
Theorem \ref{main}. The arithmetic part is the object of the paper \cite{HG}%
. The reader only interested in the $GL_{2}$ case can read a self-contained
proof of Theorem \ref{main} (both arithmetic and geometric parts) and its
consequences in this special case in \cite{BGL2}.

The novelty of the above statement resides precisely in the fact that the
integer $N(d)$ can be taken to depend only on $d$ and not on $F$ nor $%
\left\langle F\right\rangle .$ As such Theorem \ref{main} is a statement of
a different nature. What it really asserts is an inclusion of countably many
algebraic varieties into another algebraic variety. Indeed, the condition on
a $k$-tuple of matrices in say $GL_{d}(\Bbb{C})$ that they generate a
virtually solvable group is an algebraic one (see Prop. \ref{algsol} below).
On the other hand to say that no two words of length at most $N(d)$ with
letters in this $k$-tuple are generators of a free group is itself a
countable union of algebraic conditions. This way of interpreting the result
allows to derive, via an effective Nullstellensatz, several corollaries
about the girth in finite simple groups of Lie type, as well as some
diophantine properties of non-discrete subgroups of $GL_{n}(\Bbb{C)}$, in
the spirit of the works of Kaloshin-Rodnianski, Helfgott and
Bourgain-Gamburd (\cite{KR}, \cite{Hel}, \cite{BG1}, \cite{BG2}).\newline

\textit{Comments on the proof.}

Tits' proof of his alternative consists of two parts. In a first arithmetic
step, he exhibits a semisimple element of $\left\langle F\right\rangle $
which has some eigenvalue of absolute value $|\lambda |>1$ for a clever
choice of absolute value on $K$. Then in a second geometric step, he studies
the action of $\left\langle F\right\rangle $ on the projective space $\Bbb{P}%
(k^{n})$ under some suitably chosen linear representation, where $k$ is the
completion of $K$ with respect to that absolute value. The free group is
then obtained by building a so-called ping-pong pair acting on $\Bbb{P}%
(k^{n})$ (see \cite{Tits}).

The proof of Theorem \ref{main} consists in reproducing Tits' proof almost
word by word while making sure that each step can be done in a uniform way.
The arithmetic step is much harder to perform, as we need a uniform gap $%
|\lambda |>1+\varepsilon ,$ where $\varepsilon $ is allowed to depend on $d$
only. This first arithmetic step is the content of the paper \cite{HG},
which shows a height gap theorem for non amenable linear groups (see Theorem 
\ref{HeGa}). The key idea there and also in the present paper is to
introduce arithmetic heights in order to treat all absolute values of $K$ on
an equal footing. This first arithmetic step is needed only in
characteristic zero. In a second arithmetic step, we find an absolute value
for which the geometric conditions needed for the ping-pong to work are
fulfilled. This is done in Section \ref{finalproof} by estimating the
Arakelov heights of the characteristic subspaces of the matrices in $F$ in
terms of the normalized height $\widehat{h}(F)$ introduced in \cite{HG} and
by making use of another result from \cite{HG} which says that $\widehat{h}%
(F)$ can be realized up to a multiplicative factor as the height of some
conjugate of $F$ inside $SL_{d}(\overline{\Bbb{Q}})$. Once the right
absolute value has been found, the actual geometric construction of the
ping-pong pair follows Tits' geometric step very closely (unlike the
argument in \cite{uti}) ; the only notable difference is that our estimates
need to be uniform over all local fields. This requires a bit of care and is
performed in Sections \ref{proximality} and \ref{PP}. \newline

\textit{Some consequences.}

Theorem \ref{main} admits several consequences about the structure of
non-amenable linear groups. The first is a gap for the growth exponent,
namely:

\begin{corollary}
\label{growth}(\textit{Uniform exponential growth}) For every $d\in \Bbb{N}$%
, there exists a constant $\varepsilon =\varepsilon (d)>0$ such that if $K$
is any field and $F$ is a finite subset of $GL_{d}(\Bbb{C})$ containing $1$
and generating a non amenable subgroup, then for all $n\geq 1$%
\begin{equation*}
|F^{n}|\geq (1+\varepsilon )^{n}
\end{equation*}
Hence 
\begin{equation*}
\rho _{F}=\lim_{n\rightarrow +\infty }\frac{1}{n}\log |F^{n}|\geq \log
(1+\varepsilon )>0
\end{equation*}
\end{corollary}

\begin{remark}
It is possible that the assumption ``non-amenable'' in the above corollary
can be replaced by ``of exponential growth''. However we observed in \cite
{Bre} that this would imply the Lehmer conjecture about the Malher measure
of algebraic numbers. We also observed there that although every linear
solvable group of exponential growth contains a free semigroup, no analog of
Theorem \ref{main} holds for solvable groups, namely one may find sets $F_{n}
$ in $GL_{2}(\Bbb{C})$ containing $1$ and generating a solvable subgroup of
exponential growth, such that no pair of elements in $(F_{n})^{n}$ may
generate a free semigroup.
\end{remark}

\begin{remark}
\label{GrH}Examples due to Grigorchuk and de la Harpe \cite{GH} (see also 
\cite{BarCor}) show that there is a sequence of groups $\Gamma _{n}$ with
finite generating set $F_{n}$ which are virtually a direct product of
finitely many copies of the free group $F_{2}$ such that $\rho
_{F_{n}}\rightarrow 0$ as $n\rightarrow +\infty .$ Those examples can be
embedded in $SL_{m}(\Bbb{Z})$ for some possibly large $m=m(n).$ Therefore we
must have $N(d)\rightarrow +\infty $ and $\varepsilon (d)\rightarrow 0$ as $%
d\rightarrow +\infty $ in Theorem \ref{main} and Corollary \ref{growth}.
\end{remark}

The following corollary says that non-amenable linear groups have few
relations:\ there is a co-growth gap.

\begin{corollary}
\label{cogrowth}(\textit{Co-growth gap}) For every $d,k\in \Bbb{N}$, there
is $\varepsilon >0$ such that if $K$ is a field and $F=\{a_{1},...,a_{k}\}$
generates a non virtually solvable subgroup of $GL_{d}(K),$ then for every $%
n\in \Bbb{N}$, the proportion of relations $w(a_{1},...,a_{k})=1$ in the
free group $F_{k}$ of word length at most $n$ among all elements in $F_{k}$
of word length at most $n$ is at most $\exp (-\varepsilon n).$
\end{corollary}

Von Neumann showed that groups containing a free subgroup are non amenable,
i.e. have a spectral gap in $\Bbb{\ell }^{2}.$ The uniformity in Theorem \ref
{main} implies also a uniformity for the spectral gap (see \cite{Sh} for
this observation). More precisely:

\begin{corollary}
\label{gap}(\textit{Uniform Spectral Gap in }$\ell ^{2}$) For every $d\in 
\Bbb{N}$, there is $\varepsilon =\varepsilon (d)>0$ with the following
property. If $K$ is a field and $F$ is a finite subset of $GL_{d}(K)$
containing the identity and generating a non amenable subgroup and if $%
\Gamma $ is any countable subgroup of $GL_{d}(K\Bbb{)}$ containing $F$ and $%
f\in \ell ^{2}(\Gamma )$, then there is $\sigma \in F$ such that 
\begin{equation*}
\sum_{x\in \Gamma }\left| f(\sigma ^{-1}x)-f(x)\right| ^{2}\geq \varepsilon
\cdot \sum_{x\in \Gamma }|f(x)|^{2}
\end{equation*}
In particular, if $F$ in $GL_{d}(K)$ is a finite subset containing the
identity and generating a non amenable subgroup, then for every finite
subset $A$ in $GL_{2}(K)$, we have $|FA|\geq (1+\varepsilon )|A|.$
\end{corollary}

This shows also that if $\mu $ is a uniform probability measure on a set $F$
of cardinal $k$ in $GL_{d}(K)$, then the \textit{Kesten spectral radius} of $%
\mu $ (see \cite{Kes}) is uniformly bounded away from $1$ by a bound
depending only on $k$ and $d$. Hence the return probability of the simple
random walk on the group $\left\langle F\right\rangle $ decays exponentially
with an exponential rate depending only on $k$ and $d.$

The uniformity in Theorem \ref{main} allows to reduce mod $p$ and we obtain
a statement giving a lower bound on the girth of subgroups of $GL_{d}$ in
positive characteristic:

\begin{corollary}
\label{poschar}(\textit{Large girth}) Given $d,k\in \Bbb{N}$, there is $%
N_{0},N\in \Bbb{N}$ and $\varepsilon _{0},C>0$ such that for every prime $p$
and every field $K$ of characteristic $p$ and any finite $k$-element subset $%
F$ generating a subgroup of $GL_{d}(K)$ which contains no solvable subgroup
of index at most $N,$ then $F^{N_{0}}$ contains two elements $a,b$ such that 
$w(a,b)\neq 1$ in $GL_{d}(K)$ for any non trivial word $w$ in $F_{2}$ of
length at most $f(p)=C\cdot (\log p)^{\varepsilon _{0}}.$
\end{corollary}

\begin{corollary}
\label{exp}(\textit{Expansion of small sets}) There is $\varepsilon
=\varepsilon (d)>0$ such that given $k,N\in \Bbb{N}$, there is a constant $%
C_{k,N,d}$ such that for any field $K$ of charateristic $p>1$ and any subset 
$F$ of $GL_{d}(K)$ with $k$ elements generating a subgroup which has no
solvable subgroup of index at most $N$, we have $\max_{f\in
F}|A\bigtriangleup fA|\geq \varepsilon |A|$ for all subsets $A$ in $GL_{d}(K)
$ with $|A|\leq C_{k,N,d}\log \log \log p.$ 
\end{corollary}

It was conjectured in \cite{Gam} that the statement of Corollary \ref
{poschar} holds for generating subsets $F$ of $GL_{2}(\Bbb{F}_{p})$ with $%
\varepsilon _{0}=1.$ It was also proved there that a random $k$-regular
Cayley graph of $GL_{2}(\Bbb{F}_{p})$ has girth at least $(1-o(1))\log
_{k-1}(p)$.

In a similar fashion one can derive the following weak diophantine property
for subgroups of $GL_{d}(\Bbb{C})$. Let $d$ be some Riemannian distance on $%
GL_{d}(\Bbb{C}).$

\begin{corollary}
\label{wdioph}\textit{(Weak diophantine condition)} Given $d\in \Bbb{N}$,
there is $N_{0}\in \Bbb{N}$ and $\varepsilon _{1}>0$ with the following
property. For every finite set $F\subset GL_{d}(\Bbb{C})$ generating a non
virtually solvable subgroup, there is $\delta _{0}(F)>0$ such that for every 
$\delta \in (0,\delta _{0})$ there are two short words $a,b\in F^{N_{0}}$
such that $d(w(a,b),1)\geq \delta $ for every reduced word $w$ in the free
group $F_{2}$ with length $\ell (w)$ at most $(\log \delta
^{-1})^{\varepsilon _{1}}.$
\end{corollary}

In \cite{KR} Kaloshin and Rodnianski proved that for $G=SU(2,\Bbb{R)}\leq
GL_{2}(\Bbb{C})$ almost every pair $(a,b)\in G\times G$ satisfies $%
d(w(a,b),1)\geq \exp (-C(a,b)\cdot \ell (w)^{2})$ for all $w\in
F_{2}\backslash \{e\}$ and some constant $C(a,b)>0.$ Besides it is easy to
see that if $a,b\in GL_{2}(\overline{\Bbb{Q}})$ then the pair $(a,b)$
satisfies the stronger diophantine condition $d(w(a,b),1)\geq \exp
(-C(a,b)\cdot \ell (w)).$ It is conjectured in \cite{Sar} and \cite{Gam},
that this stronger condition also holds for almost every pair $(a,b)\in SU(2,%
\Bbb{R}).$

Our result also allows us to estimate the number of words of length $\leq n$
that fall in a shrinking neighborhood of $1$ in $GL_{d}(\Bbb{C}).$ More
precisely,

\begin{corollary}
\label{weq}\textit{(Weak equidistribution) }Given $d\in \Bbb{N}$, there are $%
\tau ,\varepsilon _{1},C>0$ with the following property.\textit{\ }For every 
$\{a,b\}\leq GL_{d}(\Bbb{C})$ which generates a non virtually solvable
subgroup, there is $\delta _{0}(a,b)>0$ such that for every $\delta \in
(0,\delta _{0})$ and every $n\leq C(\log \delta ^{-1})^{\varepsilon _{1}}$,
the proportion of elements $w$ in the free group $F_{2}$ of word length $n$
such that $d(w(a,b),1)\leq \delta $ is at most $\exp (-\tau n).$
\end{corollary}

In \cite{GJS}, Gamburd, Jacobson and Sarnak, showed for $G=SU(2,\Bbb{R})$
that if a pair $(a,b)\in G$ satisfies the conclusion of Corollary \ref{weq}
with $\varepsilon _{1}=1$ and $C>C_{0}$ (for some explicit $C_{0}>0$) then $%
(a,b)$ has a spectral gap on $\Bbb{L}^{2}(G).$ In \cite{BG1}, Bourgain and
Gamburd showed that if a pair $(a,b)\in G$ satisfies the above condition
with $\varepsilon _{1}=1$ and some $C=C(a,b)>0,$ then $(a,b)$ has a spectral
gap on $\Bbb{L}^{2}(G).$ This latter condition is automatically satisfied if 
$(a,b)$ satisfies the stronger diophantine condition above, for instance if $%
(a,b)\in GL_{2}(\overline{\Bbb{Q}}).$\newline

\begin{remark}
Corollaries \ref{poschar}, \ref{wdioph} and \ref{weq}, are derived from
Theorem \ref{main} by using we use a standard version of the effective
Nullstellensatz due to Masser and Wustholz (see \cite{MW}) after
reformulating Theorem \ref{main} in terms of inclusion of algebraic
varieties. See Section \ref{Corrrr}.\newline
\end{remark}

\setcounter{tocdepth}{1}
\tableofcontents
%

\section{Minimal norm and spectral radius formula}

In this section, we recall results obtained in \cite{HG} about the spectral
radius of a finite set of matrices. Given a local field $k$, we defined the
standard norm $||\cdot ||_{k}$ on $k^{d}$ to be the canonical Euclidean
(resp. Hermitian) norm if $k$ is $\Bbb{R}$ (resp. $\Bbb{C}$) and the sup
norm if $k$ is ultrametric. This induces on operator norm on the space of $%
d\times d$ matrices $M_{d}(k)$, which we again denote by $||\cdot ||_{k}.$
Given a finite subset $F$ of matrices in $M_{d}(k)$, we define its norm $%
||F||_{k}$ to be the maximal norm of any given element of $F.$ We define the
following quantities

\begin{equation*}
E_{k}(F)=\inf_{g\in GL_{d}(\overline{k})}||gFg^{-1}||_{k}
\end{equation*}
\begin{equation*}
\Lambda _{k}(F)=\max \{|\lambda |_{k},\lambda ~\text{eigenvalue}~\text{of}~%
\text{some}~f\in F\}
\end{equation*}
where $\overline{k}$ is an algebraic closure of $k$ and $|\cdot |_{k}$ is
the absolute value on $k$ extended (uniquely) to $\overline{k}.$ We also set
the \textit{spectral radius} of $F$ to be: 
\begin{equation*}
R_{k}(F)=\lim_{n\rightarrow +\infty }||F^{n}||_{k}^{\frac{1}{n}}
\end{equation*}
These quantities enjoy the following key properties.

\begin{lemma}
\label{NvsS}(Spectral Radius Formula for $F,$ \cite{HG}, Lemma 2.1.)

(a) If $k$ is ultrametric, then for any compact set $F$ containing $1$ in $%
M_{d}(k),$ there is a positive integer $q\leq d^{2}$ such that $\Lambda
_{k}(F^{q})=E_{k}(F)^{q}.$ In particular, $E_{k}(F)=R_{k}(F)=\max_{1\leq
q\leq d^{2}}\Lambda _{k}(F^{q})^{\frac{1}{q}}.$

(b) If $k$ is archimedean, there is a constant $c=c(d)\in (0,1)$ such that
for any compact set $F$ in $M_{d}(k),$ there is a positive integer $q\leq
d^{2}$ such that $\Lambda _{k}(F^{q})\geq c\cdot E_{k}(F)^{q}$. In
particular, $c\cdot E_{k}(F)\leq \max_{1\leq q\leq d^{2}}\Lambda
_{k}(F^{q})^{\frac{1}{q}}\leq R_{k}(F)\leq E_{k}(F).$
\end{lemma}

\begin{remark}
This lemma expresses in a condensed form some ideas present in the proof of
the main result of \cite{EMO} by Eskin-Mozes-Oh. It is useful to produce
elements with large eigenvalues in $F^{n}$ for some small $n$.
\end{remark}

We also record the following:

\begin{lemma}
\label{gr}(\cite{HG}, Proposition 2.5.) Suppose $k$ is archimedean (i.e. $k$
is $\Bbb{R}$ or $\Bbb{C}$). Then for every $n\in \Bbb{N}$ and every compact
subset $F$ in $SL_{d}(k)$ containing $1$, we have 
\begin{equation*}
E_{k}(F^{n})\geq E_{k}(F)^{\sqrt{\frac{n}{8d}}}
\end{equation*}
\end{lemma}

\section{Normalized height and Height gap\label{normhei}}

In this section we recall results obtained in \cite{HG} about heights. In 
\cite{HG}, we introduced the notion of normalized height $\widehat{h}(F)$ of
a finite subset of matrices $F$ in $SL_{d}(\overline{\Bbb{Q}}).$ A similar
definition can be made over the algebraic closure $\overline{\Bbb{F}_{p}(t)}$
of $\Bbb{F}_{p}(t).$ Below we recall the relevant definitions and notations.

Let $\Omega $ be either $\overline{\Bbb{Q}}$ or $\overline{\Bbb{F}_{p}(t)}$
for some prime $p>1.$ By a \textit{global field} $K,$ we mean a field
isomorphic to a finite algebraic extension of $K_{0}$, where either $K_{0}=%
\Bbb{Q}$ or $K_{0}=\Bbb{F}_{p}(t)$ for some prime $p>1.$ We denote by $V_{K}$
the set of equivalence classes of non trivial absolute values on $K$. We
make the following standard choice of representatives $|\cdot |_{v}$ for
each $v$ in $V_{K}$. Every $v\in V_{K}$ induces on $K_{0}$ an absolute value 
$v_{0}\in V_{K_{0}}.$ We first determine a set of representatives of $%
V_{K_{0}},$ then pick in each $v\in V_{K}$ the representative with that
normalization. If $K_{0}=\Bbb{Q}$ then any absolute value is equivalent to
either the standard absolute value over $\Bbb{R}$ or the $p$-adic absolute
value normalized so that $|p|_{p}=\frac{1}{p}.$ These form our
representatives. If $K_{0}=\Bbb{F}_{p}(t),$ then every absolute value is
equivalent to either $|\frac{P}{Q}|_{0}=p^{\deg P-\deg Q}$ or $|\frac{P}{Q}%
|_{\pi }=p^{(v_{\pi }(Q)-v_{\pi }(P))\cdot \deg \pi }$ where $\pi \in \Bbb{F}%
_{p}[t]$ is a monic irreducible polynomial and $v_{\pi }(P)$ is the
valuation of $\pi $ in the prime factor decomposition of $P\in \Bbb{F}%
_{p}[t].$ These form our representatives. For background on these issues see
Weil's book \cite{Weil}.

Each $v\in V_{K}$ gives rise to a local field $K_{v}$ which is the
completion of $K$ according to this absolute value. Let $n_{v}$ be the
dimension of $K_{v}$ over the closure of $K_{0}$ in $K_{v}.$ The product
formula reads 
\begin{equation}
\sum_{v\in V_{K}}n_{v}\log |x|_{v}=0  \label{pf}
\end{equation}
for every $x\in K.$ We can now recall the definition of the standard Weil
height of an algebraic number. Let $x\in K\backslash \{0\},$%
\begin{equation*}
h(x)=\frac{1}{[K:K_{0}]}\sum_{v\in V_{K}}n_{v}\log ^{+}|x|_{v},
\end{equation*}
where $\log ^{+}=\max \{\log ,0\}.$

In \cite{HG}, we introduced the following heights for $F$ a finite subset of 
$M_{d}(K)\backslash \{0\},$%
\begin{equation}
h(F)=\frac{1}{[K:K_{0}]}\sum_{v\in V_{K}}n_{v}\log ^{+}||F||_{v}
\label{defheight}
\end{equation}
We also defined the \textit{normalized height of }$F$ as 
\begin{equation*}
\widehat{h}(F)=\lim_{n\rightarrow +\infty }\frac{1}{n}h(F^{n})=\frac{1}{%
[K:K_{0}]}\sum_{v\in V_{K}}n_{v}\log ^{+}R_{v}(F)
\end{equation*}
and the \textit{minimal height of }$F$ as 
\begin{equation*}
e(F)=\frac{1}{[K:K_{0}]}\sum_{v\in V_{K}}n_{v}\log ^{+}E_{v}(F)
\end{equation*}
where we have denoted by $E_{v}(F)$ (resp. $R_{v}(F)$) the quantity $%
E_{K_{v}}(F)$ (resp. $R_{K_{v}}(F)$) defined above. Observe that the height $%
h(F)$ depends on the choice of basis in $K^{d},$ while the normalized height 
$\widehat{h}(F)$ and minimal height $e(F)$ do not. We will often write $%
h=h_{f}+h_{\infty }$ to distinguish the finite part and the infinite part of
the height in the obvious way.

In \cite{HG} we proved the following results:

\begin{lemma}
(\cite{HG} Proposition 2.18) There is a constant $c_{1}=c_{1}(d)>0$ such
that for every finite subset $F$ in $M_{d}(\Omega )$%
\begin{equation*}
e(F)\geq \widehat{h}(F)\geq e_{f}(F)+c_{1}\cdot e_{\infty }(F)\min
\{1,e_{\infty }(F)\}
\end{equation*}
\end{lemma}

It is easy to verify that $\widehat{h}(F)=0$ if and only if $e(F)=0$ if and
only if $\left\langle F\right\rangle $ is virtually unipotent. Note in
particular that if $char(\Omega )>0$, there are no infinite places so the
normalized height and the minimal height coincide, and

\begin{lemma}
\label{zerocase}If $char(\Omega )=p>0$ then $\widehat{h}(F)=0$ if and only
if $\left\langle F\right\rangle $ finite.
\end{lemma}

\proof%
The if part follows from the definition of $\widehat{h}(F).$ Suppose now
that $\widehat{h}(F)=0$. Then for every eigenvalue $\lambda $ of an element $%
g\in F,$ $h(\lambda )=0,$ hence $\lambda $ is of finite order and belongs to 
$\overline{\Bbb{F}_{p}}$. Hence $g$ also is of finite order as both the
semisimple part $g_{s}$ and the unipotent part $g_{u}$ are of finite order.
But Shur's theorem (see \cite{CurR}) says that any finitely generated
torsion linear group is finite. 
\endproof%

The main theorem of \cite{HG} is the following.

\begin{theorem}
\label{HeGa}(Height Gap, \cite{HG} Theorem 1.1) There is a constant $%
\varepsilon =\varepsilon (d)>0,$ such that if $F$ is a finite subset of $%
SL_{d}(\overline{\Bbb{Q}})$ generating a non virtually solvable subgroup,
then 
\begin{equation*}
\widehat{h}(F)>\varepsilon 
\end{equation*}
\end{theorem}

Given a Chevalley group $\Bbb{G}$, there is a special choice of basis of the
Lie algebra $\frak{g}=Lie(\Bbb{G})$ which is made of weight vectors of a
maximal split torus and defines a $\Bbb{Z}$-structure on $\Bbb{G}$ (see
Steinberg's notes \cite{Stein}, and Paragraph \ref{irredrep} below). With
respect to this basis and viewing $\Bbb{G}$ as a subgroup of $SL_{d}(\frak{g}%
)$ we may define the height $h(g)$ for any $g\in \Bbb{G}(\Omega )$ as in $(%
\ref{defheight})$. We then have:

\begin{theorem}
\label{GoodPos}(\cite{HG} Proposition 3.3) If $\Bbb{G}$ is a Chevalley
group, then there is a constant $C=C(\Bbb{G})>0$ and a Zariski open subset $%
\mathcal{O}=\mathcal{O}(\Bbb{G})$ of $\Bbb{G\times G}$ such that for any
choice of $\Omega $ and for any pair $(a,b)\in \mathcal{O}(\Omega ),$ there
is $g\in \Bbb{G}(\Omega )$ such that, setting $F=\{a,b\}$, 
\begin{equation}
h(gFg^{-1})\leq C\cdot \widehat{h}(F)  \label{RPos}
\end{equation}
\end{theorem}

Unlike Theorem \ref{GoodPos}, there is no analog of Theorem \ref{HeGa} for $%
\Omega =\overline{\Bbb{F}_{p}(t)}$. In \cite{HG}, we proved Theorem \ref
{GoodPos} when $\Omega =\overline{\Bbb{Q}}$ because we were only concerned
with characteristic zero. However the proof we gave works the same word by
word in the positive characteristic case, and is even simpler since in that
case there are no infinite places : in particular the additive constants $%
C_{\infty }$ and $C_{\infty }^{\prime }$ that are obtained along the way
vanish and the use of Theorem \ref{HeGa} to get rid of them is not needed
(see \cite{HG}).

\subsection{Arakelov Height on Grassmannians\label{Arak}}

Here we record some well-known facts about Arakelov heights. Let $K$ be a
global field. The Arakelov height on the projective space $\Bbb{P}(K^{d})$
is defined as follows (see \cite{Bom}) for $x=(x_{1}:...:x_{d})$, 
\begin{equation*}
h_{Ar}(x)=\frac{1}{[K:K_{0}]}\sum_{v\in V_{K}}n_{v}\log ||x||_{v}
\end{equation*}
where $||x||_{v}$ is the standard norm on $K_{v}^{d}$ as defined above. It
is well defined thanks to the product formula (\ref{pf}) and always
non-negative. This allows to define the height of a projective linear
subspace of $\Bbb{P}(K^{d}).$ Indeed if $W\leq \Bbb{P}(K^{d})$ is such then
we set 
\begin{equation*}
h_{Ar}(W)=h_{Ar}(\Lambda ^{\dim W}W)
\end{equation*}
where $\Lambda ^{\dim W}W$ is the wedge product of $W$ viewed as a
projective point in the projective space $\Bbb{P}(\Lambda ^{\dim W}K^{d})$.
By convention we set $h_{Ar}(\{0\})=0$. Recall that the following holds for
two projective linear subspaces (see \cite{Bom}) $V$ and $W$%
\begin{equation*}
h_{Ar}(V)+h_{Ar}(W)\geq h_{Ar}(V+W)+h_{Ar}(V\cap W).
\end{equation*}
Moreover for every linear form $f$, seen as a point in the dual space $%
(K^{d})^{*}$, $h_{Ar}(f)$ makes sense as we have (see \cite{Bom}), 
\begin{equation*}
h_{Ar}(\ker f)=h_{Ar}(f)
\end{equation*}
and more generally, $h_{Ar}(W)=h_{Ar}(W^{\bot })$, where $W^{\bot }$ is the
orthogonal of $W$ in $(K^{d})^{*}.$

Also note that if $g\in SL_{d}(K)$ and $W$ is a subspace of $K^{d},$ then 
\begin{equation*}
h_{Ar}(gW)\leq d\cdot h(g)+h_{Ar}(W)
\end{equation*}
where $h(g)=h(\{g\})$ as defined in the last paragraph. Note also that $%
h(g^{-1})\leq (d-1)h(g).$

\begin{definition}
\label{admiss}Given $A\in SL_{d}(K),$ we will say that a vector subspace $W$
(or its projectivization) is $A$\textbf{-admissible} if it is a sum of
generalized eigenspaces of $A.$ We also denote by $W^{c}$ its complementary
subspace, i.e. the sum of the remaining generalized eigenspaces.
\end{definition}

\begin{lemma}
\label{heightadmis}Let $W$ be an $A$-admissible subspace. Then 
\begin{equation*}
h_{Ar}(W)\leq d^{2}\cdot (2h(A)+\varepsilon _{\Omega }\log 2)
\end{equation*}
where $\varepsilon _{\Omega }=0$ if $car(K)>1$ and $1$ if $car(K)=0.$
\end{lemma}

Note that if $A,B\in M_{d}(K),$ then $h(A+B)\leq h(A)+h(B)+\varepsilon
_{\Omega }\log 2.$ Moreover, if $\alpha $ is an eigenvalue of $A$, then $%
h(\alpha )\leq h(A).$ Also $h_{Ar}(\ker A)\leq (rk(A))\cdot h(A).$ Indeed $%
h_{Ar}(\ker A)=h_{Ar}((\ker A)^{\bot })=h_{Ar}(\func{Im}A^{t}).$ But for $%
B\in M_{d}(K),$ $h_{Ar}(\func{Im}B)\leq (rk(B))\cdot h(B),$ since we may
choose a subset of the canonical basis, say $e_{1},...,e_{k}$ such that $%
Be_{1},...,Be_{k}$ generates $\func{Im}B,$ and $||Be_{1}\wedge ...\wedge
Be_{k}||_{v}\leq ||B||_{v}^{k}$ where $k=rk(B).$ With these observations in
hand we can prove Lemma \ref{heightadmis}.

\proof%
We have $W=\bigoplus E_{\alpha }$ for some eigenvalues $\alpha $ of $A,$
where $E_{\alpha }$ is the corresponding generalized eigenspace. Hence $%
h_{Ar}(W)\leq \sum h_{Ar}(E_{\alpha }).$ If $n_{\alpha }=\dim E_{\alpha },$
then $E_{\alpha }=\ker (A-\alpha )^{n_{\alpha }}.$ Hence $h_{Ar}(E_{\alpha
})\leq d\cdot h((A-\alpha )^{n_{\alpha }})=dn_{\alpha }\cdot
(2h(A)+\varepsilon _{\Omega }\log 2).$ Hence the result. 
\endproof%

\section{Proximality\label{proximality}}

In this paragraph we recall the well-known notion of a proximal element in $%
SL_{d}(k),$ where $k$ is a local field, and we show some precise estimates
as to how such elements act on the projective space $\Bbb{P}(k^{d}).$ The
results of this paragraph are contained in Lemma \ref{prox} and Lemma \ref
{TitsConv} below.

A element $a\in SL_{d}(k)$ is said to be \textbf{proximal} if there is a
unique (multiplicity one) eigenvalue of $a$ with maximum modulus $\Lambda
_{k}(a).$ We will also need to consider almost proximal elements where the
eigenvalues which are larger than, say, some $\omega $ are much larger than
all other eigenvalues.

Lemma \ref{prox} computes the rate of convergence to the attracting point of
powers of a given proximal element $a$ in terms of three quantities : its
norm $||a||,$ the modulus of its maximal eigenvalue $\Lambda _{k}(a)$ and
the modulus of its second to maximal eigenvalue $\lambda _{k}(a).$ A similar
estimate is given for an almost proximal element depending on the choice of
the cursor $\omega .$ Lemma \ref{TitsConv} is a converse statement
originally used by Tits in the proof of his alternative which gives a
sufficient condition for $a\in SL_{d}(k)$ to be proximal : it is as soon as $%
a$ stabilizes some open subset where it contracts distances.

We had to be careful in those estimates, and they differ in some non
insignificant ways from the estimates used in earlier works (as in \cite{AMS}%
\cite{uti}). In particular they are uniform over all ultrametric local
fields. The multiplicative constants $C_{k,i}$'s that appears in the
estimates always disappears when $k$ is ultrametric. This will turn out to
be crucial for us in the sequel.

\subsection{The Fubini-Study metric on $\Bbb{P}(k^{d})\label{FubStudSec}$}

Let $k$ be a local field and $\overline{k}$ an algebraic closure of $k$.
Recall that we endow the projective space $\Bbb{P}(\overline{k}^{d})$ with
the standard (Fubini-Study) distance defined by 
\begin{equation}
d([u],[v])=\frac{||u\wedge v||}{||u||\cdot ||v||}  \label{FubStud}
\end{equation}
for any $u,v\in \overline{k}^{d}\backslash \{0\}$ and $||\cdot ||$ is the
standard norm on $\overline{k}^{d}$ (i.e. Euclidean norm if $k$ is
archimedean and sup norm if $k$ is non archimedean). To avoid heavy
notation, we will denote by the same letter a non zero vector, or subspace
of $\overline{k}^{d}$ and its projectivization in $\Bbb{P}(\overline{k}%
^{d}). $ This ambiguity should not lead to any serious confusion.

We denote by $\Bbb{K}_{k}$ the maximal compact subgroup of $SL_{d}(k)$ equal
to $SO(d,\Bbb{R})$ if $k=\Bbb{R}$, $SU(d,\Bbb{R})$ if $k=\Bbb{C}$, and $%
SL_{d}(\mathcal{O}_{k})$ is $k$ is ultrametric. Its action on $\Bbb{P}%
(k^{d}) $ preserves $d$ (in fact this characterizes $d$ up to composition by
some positive function).

\begin{lemma}
\label{lip}Let $h\in SL_{d}(k).$ Then $Lip(h)\leq (||h||\cdot
||h^{-1}||)^{2}\leq ||h||^{2d},$ where $Lip(h)$ is the smallest constant $%
L\geq 0$ such that $d(hx,hy)\leq L\cdot d(x,y)$ for all $x,y\in \Bbb{P}%
(k^{d}).$
\end{lemma}

\proof%
Writing $h$ in Cartan's $\Bbb{K}_{k}A\Bbb{K}_{k}$ decomposition of $%
SL_{d}(k),$ one sees that we can assume that $h$ is diagonal and we are thus
reduced to a straightforward verification.%
\endproof%

Recall that if $H$ is a hyperplane in $k^{d},$ and $f$ a non zero linear
form on $k^{d}$ with kernel $H,$ then if $u\in k^{d}\backslash \{0\},$ its
distance to $H$ is 
\begin{equation}
d(u,H)=\frac{|f(u)|}{||f||\cdot ||u||}  \label{disttoH}
\end{equation}
where $||f||=\sup \{|f(x)|,||x||\leq 1,x\in k^{d}\}.$ More generally, if $V$
and $W$ are two $k$-subspaces in direct sum, i.e. $V\oplus W=k^{d},$ then 
\begin{equation}
d(V,W)=\frac{||\underline{v}\wedge \underline{w}||}{||\underline{v}||\cdot ||%
\underline{w}||}  \label{distsubspace}
\end{equation}
where $\underline{v}=v_{1}\wedge ...\wedge v_{l}$ and $\underline{w}%
=w_{1}\wedge ...\wedge w_{d-l}$ for any basis $(v_{1},...,v_{l})$ of $V$ and 
$(w_{1},...,w_{d-l})$ of $W.$ In particular, when $k$ is archimedean, two
subspaces are orthogonal if and only if they are at distance $1$. Let $%
(e_{1},...,e_{d})$ be the canonical basis in $k^{d}.$

\begin{lemma}
\label{dist0}Let $f$ be a non-zero linear form on $k^{d}$ and $H=\ker f$.
Let $V$ a $k$-subspace in $k^{d}$ and $V^{*}$ the orthogonal of $V$ in the
dual of $k^{d}.$ Then for every $v\in V,$%
\begin{equation}
d(v,H)=d(v,V\cap H)\cdot d(f,V^{*})  \label{distsub}
\end{equation}
\end{lemma}

\proof%
Observe that as $\Bbb{K}_{k}$ permutes transitively the $k$-subspaces of
given dimension and preserves $d,$ we may assume that $V=\left\langle
e_{1},...,e_{p}\right\rangle $ for some $p\in [0,d].$ Then we may write $f$
in the dual canonical basis $f=\sum f_{i}e_{i}^{*}=f^{<}+f^{>}$ where $f^{<}$
is the part of the sum involving indices $i\leq k$ and $f^{>}$ the other
part. Let $\underline{e}^{>}=e_{p+1}^{*}\wedge ...\wedge e_{d}^{*}$. Then $%
||f||\cdot d(f,V^{*})=||f\wedge \underline{e}||=||f^{<}\wedge \underline{e}%
||=||f^{<}||.$ On the other hand note that $f^{<}$ coincides with $f$ on $V.$
Hence for $v\in V$, $d(v,V\cap H)=\frac{|f(v)|}{||f^{<}||\cdot ||v||}.$ As $%
d(v,H)=\frac{|f(v)|}{||f||\cdot ||v||}$, combining these relations we do
obtain $(\ref{distsub})$. 
\endproof%

\begin{lemma}
\label{dist1}Let $V\oplus W=k^{d}$ and $H$ a hyperplane in $k^{d}$ with $%
V\nsubseteq H.$ Let $\pi $ be the linear projection onto $V$ with kernel $W$%
. Then for every $u\in \Bbb{P(}k^{d})\backslash W$ we have 
\begin{equation*}
d(\pi u,V\cap H)\geq d(u,W+V\cap H)\cdot d(V,W)
\end{equation*}
\end{lemma}

\proof%
Write $u=\pi u+\pi u^{\bot }\in V\oplus W$. If $v_{1},...,v_{k-1}$ is a
basis of $V\cap H$ and $w_{1},...,w_{d-k}$ a basis of $W$ we set $\underline{%
v}=v_{1}\wedge ...\wedge v_{k-1}$ and $\underline{w}=w_{1}\wedge ...\wedge
w_{d-k}$. We have $d(\pi u,V\cap H)\geq d(\pi u,W+V\cap H)=\frac{||\pi
u\wedge \underline{v}\wedge \underline{w}||}{||\pi u||\cdot ||\underline{v}%
\wedge \underline{w}||}=d(u,W+V\cap H)\cdot \frac{||u||}{||\pi u||}.$ We may
assume $u\notin V$, then on the other hand $d(V,W)\leq d(\pi u,\pi u^{\bot
})=\frac{||u\wedge \pi u^{\bot }||}{||\pi u||\cdot ||\pi u^{\bot }||}\leq 
\frac{||u||}{||\pi u||}.$ We are done. 
\endproof%

\subsection{Contraction properties of proximal and almost proximal elements}

For $a\in SL_{d}(k)$ we set $E_{\lambda }$ its generalized eigenspace with
eigenvalue $\lambda $. In this paragraph, we will assume that eigenvalues of 
$a$ belong to $k.$ We let $\Lambda _{k}(a)=\max \{|\mu |_{k},\mu $
eigenvalue of $a\}$ and $\lambda _{k}(a)$ the modulus of the second heighest
eigenvalue of $a$. An element $a\in SL_{d}(k)$ is said to be \textbf{proximal%
} if $\lambda _{k}(a)<\Lambda _{k}(a).$

To deal with non proximal elements we introduce some positive real number $%
\omega >0$, such that $\Lambda _{k}(a^{-1})^{-1}<\omega \leq \Lambda
_{k}(a). $ We set $\Lambda _{k}^{\omega }(a)=\min \{|\mu |_{k},\mu $
eigenvalue of $a, $ $|\mu |_{k}\geq \omega \}$ and $\lambda _{k}^{\omega
}(a)=\max \{|\mu |_{k},\mu $ eigenvalue of $a,|\mu |_{k}<\omega \}.$

\begin{lemma}
\label{almostprox}Suppose $a\in SL_{d}(k)$ and let $A=\Lambda _{k}(a)\Lambda
_{k}(a^{-1})\geq 1.$ For every $\varepsilon >0$ there is $\eta =\eta
(\varepsilon ,d)>0$ and $\omega $ such that 
\begin{equation*}
\Lambda _{k}(a^{-1})^{-1}<\omega \leq \Lambda _{k}(a)
\end{equation*}
and 
\begin{equation}
A^{\eta }\cdot \left( \frac{\Lambda _{k}(a)}{\Lambda _{k}^{\omega }(a)}%
\right) ^{\frac{1}{\varepsilon }}\leq \frac{\Lambda _{k}^{\omega }(a)}{%
\lambda _{k}^{\omega }(a)}  \label{almostproxineq}
\end{equation}
\end{lemma}

\proof%
Let $\lambda _{1},...,\lambda _{d}$ be the eigenvalues of $a$ ordered as $%
|\lambda _{1}|_{k}\geq ...\geq |\lambda _{d}|_{k}.$ Let $\ell _{i}=\log 
\frac{|\lambda _{i}|}{|\lambda _{i+1}|}\geq 0.$ Fix $\varepsilon >0$ and
take some $\eta >0.$ We claim that for $\eta $ small enough, there exists $%
i_{0}\in [0,d-2]$ such that $\ell _{i_{0}+1}-\eta \log A\geq \frac{1}{%
\varepsilon }\cdot (\ell _{1}+...+\ell _{i_{0}}).$ Indeed, otherwise we
would have $\ell _{1}<\eta \log A$, $\ell _{2}<\eta \log A+\frac{1}{%
\varepsilon }\ell _{1}$, etc, until we get $\log A=\ell _{1}+...+\ell
_{d-1}\leq C(\varepsilon ,d)\eta \log A$ for some computable constant $%
C(\varepsilon ,d),$ a contradiction if $\eta $ is smaller than say $\frac{1}{%
2C(\varepsilon ,d)}$. Let $\omega =|\lambda _{i_{0}+1}|_{k}.$ We are done. 
\endproof%

When $\varepsilon $ is fixed and $\omega $ so given by Lemma \ref{almostprox}%
, we will refer to $a$ as being \textit{almost proximal for }$\omega $%
\textit{.}

We will let $H_{a}^{\omega }$ be the vector subspace equal to the sum of the 
$E_{\lambda }$'s for which $|\lambda |_{k}\leq \lambda _{k}^{\omega }(a).$
Similarly, we denote its complementary subspace by $V_{a}^{\omega
}=\bigoplus E_{\lambda }$, the sum being over those $\lambda $'s such that $%
|\lambda |_{k}\geq \Lambda _{k}^{\omega }(a)$. We let $\pi _{a}^{\omega }$
be the linear projection onto $V_{a}^{\omega }$ with kernel $H_{a}^{\omega
}. $ We also set $l_{\omega }=\dim V_{a}^{\omega }.$ If $a$ is proximal, we
will drop the superscript $\omega $ (and set it to be $\Lambda _{k}(a)$) and
simply denote by $V_{a}$, $H_{a}$, and $\pi _{a}$ the corresponding
quantities.

\begin{remark}
Note that if $a\in GL_{d}(k),$ then its eigenvalues belong to the extension
of $k$ generated by all algebraic extensions of $k$ in a given algebraic
closure $\overline{k}$ of degree at most $d$ (there are finitely many such).
So this extension is also a local field. Hence up to passing to this finite
extension one may always assume that the eigenvalues of $a$ belong to $k.$
\end{remark}

Lemma \ref{prox} below is the main result of this section. Its proof will
occupy the subsequent two paragraphs. When $k$ is archimedean, let $C_{k}=2$
and $C_{k,1}=d.$ When $k$ is ultrametric set $C_{k}=C_{k,1}=1.$ Let also $%
p(d)=10^{d}$ (these are only given as crude estimates, we made no attempt at
finding sharp constants in this statement). Finally let $L_{k}^{\omega
}(a)=1 $ if $k$ is ultrametric while $L_{k}^{\omega }(a)=\frac{\Lambda
_{k}^{\omega }(a)}{\Lambda _{k}^{\omega }(a)-\lambda _{k}^{\omega }(a)}$ if $%
k$ is archimedean.

\begin{lemma}
\label{prox}Let $a\in SL_{d}(k)$ whose eigenvalues belong to $k$ and assume $%
\omega $ is a real number such that $\Lambda _{k}(a^{-1})^{-1}<\omega \leq
\Lambda _{k}(a).$ Let $l_{\omega }=\dim V_{a}^{\omega }$ and $\pi
_{a}^{\omega }$ the projection on $V_{a}^{\omega }$ with kernel $%
H_{a}^{\omega }.$ Then for any $u\neq v\in \Bbb{P(}k^{d}),$ and any integer $%
n\in \Bbb{N}$%
\begin{equation}
d(a^{n}u,\pi _{a}^{\omega }(a^{n}u))\cdot d(u,H_{a}^{\omega })\leq
(C_{k}\cdot ||a||_{k})^{p(d)}\cdot \left( C_{k,1}^{l_{\omega }}\cdot \left( 
\frac{\Lambda _{k}(a)}{\Lambda _{k}^{\omega }(a)}\right) ^{l_{\omega
}-1}\cdot \frac{\lambda _{k}^{\omega }(a)}{\Lambda _{k}^{\omega }(a)}\right)
^{n}  \label{cont}
\end{equation}
Furthermore 
\begin{equation}
\frac{d(a^{n}u,a^{n}v)\cdot d(v,H_{a}^{\omega })\cdot d(u,H_{a}^{\omega })}{%
d(u,v)}\leq (C_{k}\cdot L_{k}^{\omega }(a)\cdot ||a||_{k})^{p(d)}\cdot
\left( C_{k,1}^{2l_{\omega }+2}\cdot \left( \frac{\Lambda _{k}(a)}{\Lambda
_{k}^{\omega }(a)}\right) ^{2l_{\omega }-1}\cdot \frac{\lambda _{k}(a)}{%
\Lambda _{k}(a)}\right) ^{n}  \label{lipa}
\end{equation}
\end{lemma}

Observe that $(\ref{cont})$ says nothing if the quantity inside the bracket
is not $<1.$

The following Tits Converse Lemma is useful when one needs to build an
element $x$ such that both $x$ and $x^{-1}$ are proximal.

\begin{lemma}
\label{TitsConv}(Tits Converse Lemma \cite{Tits}) Let $a\in SL_{d}(k)$.
Assume that there exists a point $v\in \Bbb{P}(k^{d})$ and an open
neighborhood $U$ of $v$ such that $\overline{aU}\subset U$ and such that $%
Lip(a_{|U})<1$, where $Lip(a_{|U})$ is the smallest constant $L>0$ such that 
$d(ax,ay)\leq L\cdot d(x,y)$ for every $x,y\in U.$ Then $a$ is proximal, $%
V_{a}\in U$ and $\frac{\lambda _{k}(a)}{\Lambda _{k}(a)}\leq Lip(a_{|U}).$
\end{lemma}

\proof%
The compact subset $\overline{aU}$ is stable under $a$ and on it $a$
contracts distances. It follows immediately that all orbits $(a^{n}u)_{n\geq
0}$ converge to the unique fixed point $v_{a}$ of $a$ in $\overline{aU}.$
Let $\alpha $ be the eigenvalue of $a$ with eigenvector $v_{a}.$ Let $\beta $
be another eigenvalue of $a$ (if $\alpha $ has multiplicity higher than $1,$
we may take $\beta =\alpha $). There is a non zero vector $w$ such that $%
aw=\beta w+\kappa v_{a}$ for some $\kappa \in k.$ Let $\varepsilon \in
k\backslash \{0\}$ with $|\varepsilon |_{k}$ arbitrarily small. Then one
computes from (\ref{FubStud}) $\lim_{|\varepsilon |\rightarrow 0}\frac{%
d(a(v_{a}+\varepsilon w),v_{a})}{d(v_{a}+\varepsilon w,v_{a})}=\frac{|\beta |%
}{|\alpha |}.$ If $|\varepsilon |_{k}$ is small enough, $v_{a}+\varepsilon
w\in U$ and thus $\frac{|\beta |}{|\alpha |}\leq Lip(a_{|U})<1.$ We are done.%
\endproof%

\subsection{Four intermediary geometric lemmas}

In this paragraph, we state and prove four intermediary results needed in
the proof of Lemma \ref{prox}. Unless otherwise stated $a\in GL_{d}(k)$ and
its eigenvalues belong to $k$.

\begin{lemma}
\label{triang}Let $a\in GL_{d}(k)$ and $\alpha $ an eigenvalue of $a$. Then
there is some $h\in \Bbb{K}_{k}$ such that $hah^{-1}$ is a lower triangular
matrix with top left entry equal to $\alpha $.
\end{lemma}

\proof%
Since eigenvalues of $a$ belong to $k$, $a$ and hence also its transpose $%
a^{t}$ are triangularizable over $k$, i.e. $a^{t}$ stabilizes a full $k$%
-flag $\mathcal{F}$. We may also assume that $\mathcal{F}$ starts with the
line $kv,$ where $v$ is an eigenvector of $a^{t}$ with eigenvalue $\alpha .$
But full $k$-flags are conjugate under $GL_{d}(k).$ Hence $\mathcal{F}=g%
\mathcal{F}_{0}$ where $\mathcal{F}_{0}$ is the standard flag generated by
the canonical basis of $k^{d}$ and $ge_{1}=v.$ The Iwasawa decomposition
reads $GL_{d}(k)=\Bbb{K}_{k}B_{0}$ where $B_{0}$ is the Borel stabilizing $%
\mathcal{F}_{0}.$ Thus we may assume that $g\in \Bbb{K}_{k}.$ Thus $%
g^{-1}a^{t}g$ stabilizes $\mathcal{F}_{0}$ and is upper triangular. Hence $%
h=g^{t}\in \Bbb{K}_{k}$ will do. 
\endproof%

Let $C_{k,1}$ be equal to $d$ if $k$ is archimedean and equal to $1$ if $k$
is ultrametric.

\begin{lemma}
\label{Conj}Let $a\in GL_{d}(k).$ Then there exists an $h\in SL_{d}(k)$ such
that $\left\| hah^{-1}\right\| \leq C_{k,1}\Lambda _{k}(a)$ and $\max
\{\left\| h\right\| ,||h^{-1}||\}\leq \left\| a\right\| ^{\frac{d-1}{2}}.$
\end{lemma}

\proof%
By Lemma \ref{triang}, one may assume that $a\in GL_{d}(k)$ is lower
triangular. Let $h=t^{\frac{d+1}{2}}diag(t^{-1},...,t^{-d})\in SL_{d}(k).$
Choose $t\in k$ such that $|t^{-1}|_{k}=||a||_{k}.$ Then $\max \{\left\|
h\right\| ,||h^{-1}||\}\leq \left\| a\right\| ^{\frac{d-1}{2}}$ and the
off-diagonal coefficients of $hah^{-1}$ are of modulus $\leq 1.$ As $%
||a||\leq C_{k,1}\max |a_{ij}|,$ we are done. 
\endproof%

\begin{remark}
\label{Conjplane}Note that we also get $||\Lambda ^{2}(hah^{-1})||\leq
C_{k,1}^{4}\Lambda _{k}(a)\lambda _{k}(a).$
\end{remark}

Recall that $C_{k}$ is $2$ if $k$ is archimedean and $1$ if $k$ is
ultrametric.

\begin{lemma}
\label{interm} Let $a\in SL_{d}(k)$ and $\omega $ with $\Lambda
_{k}(a^{-1})^{-1}<\omega \leq \Lambda _{k}(a).$ Let $l_{\omega }=\dim
V_{a}^{\omega }$ and $L_{k}^{\omega }(a)=1$ if $k$ is ultrametric while $%
L_{k}^{\omega }(a)=\frac{\Lambda _{k}^{\omega }(a)}{\Lambda _{k}^{\omega
}(a)-\lambda _{k}^{\omega }(a)}$ if $k$ is archimedean. Then $%
d(V_{a}^{\omega },H_{a}^{\omega })^{-1}\leq (C_{k}\cdot L_{k}^{\omega
}(a)||a||^{l})^{\left( \binom{d}{l}-1\right) }.$
\end{lemma}

\proof%
First let us assume $a$ is proximal (and $\omega =\Lambda _{k}(a)$) with
eigenvalue of maximal modulus $\alpha _{1},$ and let $\alpha _{2},...,\alpha
_{d}$ be the remaining eigenvalues. By Lemma \ref{triang}, one may assume
that $a\in SL_{d}(k)$ is lower triangular with $\alpha _{1}$ in the upper
left corner. Let $V_{a}\in k^{d}\backslash \{0\}$ be such that $%
aV_{a}=\alpha _{1}V_{a}.$ As $V_{a}\notin H_{a}=\left\langle
e_{2},...,e_{d}\right\rangle ,$ we may assume that $%
V_{a}=(1,x_{2},...,x_{d}) $ in the canonical basis. Then $%
d(V_{a},H_{a})=1/||V_{a}||.$ Decomposing $aV_{a}=\alpha _{1}V_{a}$ in
coordinates, we obtain $(\alpha _{1}-\alpha _{2})v_{2}=a_{21},$ $(\alpha
_{1}-\alpha _{3})v_{3}=a_{31}+a_{32}v_{2},$ etc. This allows to recursively
estimate each $v_{i}$ and at the end we get that $||V_{a}||\leq (1+L_{k}^{2}%
\frac{||a||^{2}}{\Lambda _{k}(a)^{2}})^{(d-1)/2}\leq 2^{(d-1)/2}\cdot
L_{k}^{d-1}\cdot \left( \frac{||a||}{\Lambda _{k}(a)}\right) ^{d-1}$ when $k$
is archimedean while when $k$ is ultrametric $||V_{a}||\leq
L_{k}^{d-1}\left( \frac{||a||}{\Lambda _{k}(a)}\right) ^{d-1}.$

We now explain how to reduce the general case to the proximal case. Let $%
(v_{1},...,v_{l})$ and $(w_{1},...,w_{d-l})$ be respective basis of $%
V_{a}^{\omega }$ and $H_{a}^{\omega }.$ Let $\underline{v}=v_{1}\wedge
...\wedge v_{l}$ and $\underline{w}=w_{1}\wedge ...\wedge w_{d-l}.$ From $(%
\ref{distsubspace})$ we have $d(V_{a}^{\omega },H_{a}^{\omega })=\frac{||%
\underline{v}\wedge \underline{w}||}{||\underline{v}||\cdot ||\underline{w}||%
}.$ The canonical map $\Lambda ^{l}k^{d}\times \Lambda
^{d-l}k^{d}\rightarrow k$ establishes an isomorphism between $\Lambda
^{d-l}k^{d}$ and the dual of $\Lambda ^{l}k^{d}.$ Under this identification $%
\underline{w}$ is a linear form on $\Lambda ^{l}k^{d}$ and formulae $(\ref
{distsubspace})$ and $(\ref{disttoH})$ coincide, i.e. $d(V_{a}^{\omega
},H_{a}^{\omega })=d(\Lambda ^{l}V_{a}^{\omega },\ker \underline{w}).$ On
the other hand $\Lambda ^{l}a$ is proximal on $\Lambda ^{l}k^{d}$ with $%
V_{\Lambda ^{l}a}=\Lambda ^{l}V_{a}^{\omega }$ and $H_{\Lambda ^{l}a}=\ker 
\underline{w}.$ Hence by the above $d(\Lambda ^{l}V_{a},\ker \underline{w}%
)^{-1}\leq (C_{k}\cdot L_{k}^{\omega }(\Lambda ^{l}a)\cdot \frac{||\Lambda
^{l}a||}{\Lambda _{k}(\Lambda ^{l}a)})^{D-1}$ where $D=\dim \Lambda
^{l}k^{d}=\binom{d}{l}$ and the result follows as $\Lambda _{k}(\Lambda
^{l}a)\geq 1$. 
\endproof%

Recall that $C_{k,1}$ is $d$ if $k$ is archimedean and $1$ if $k$ is
ultrametric.

\begin{lemma}
\label{interm2}Let $a\in SL_{d}(k)$ with $\Lambda _{k}(a^{-1})^{-1}<\omega
\leq \Lambda _{k}(a)$ and $l_{\omega }=\dim V_{a}^{\omega }.$ Set $%
V=\left\langle e_{2},...,e_{l_{\omega }}\right\rangle $, $H=\left\langle
e_{l_{\omega }+1},...,e_{d}\right\rangle .$ There exists $h\in SL_{d}(k)$
with $hV_{a}^{\omega }=V,$ $hH_{a}^{\omega }=H$ such that if we set $%
\widetilde{a}=hah^{-1}$ then $\left\| \widetilde{a}_{|H}\right\| \leq
C_{k,1}\lambda _{k}^{\omega }(a)$ and $\left\| \widetilde{a}_{|V}\right\|
\leq C_{k,1}\Lambda _{k}(a)$ and $\left\| h^{-1}\right\| \leq \left( \frac{%
\sqrt{C_{k,1}}}{d(V_{a}^{\omega },H_{a}^{\omega })}\right) ^{\frac{d(d+1)}{2}%
}\left\| a\right\| ^{\frac{d-1}{2}}.$
\end{lemma}

\proof%
First note that applying Lemma \ref{triang} we can assume that $a$ is lower
triangular and that $H_{a}^{\omega }=H.$ Then observe that for any subspace $%
F$ of $k^{d},$ one may find a basis $f_{1},...,f_{p}$ of $F$ such that $%
||f_{1}\wedge ...\wedge f_{p}||=1$ and $||f_{i}||=1$ for each $i=1,...,p.$
Choose such a basis, say $v_{1},...,v_{l}$ of $V_{a}^{\omega }$ and, for $%
\mu \in k$ to be defined later, denote by $h_{1}\in GL_{d}(k)$ the map $%
h_{1}v_{i}=e_{i}$ if $i<l$, $h_{1}v_{l}=\mu e_{l}$ and $h_{1}e_{i}=e_{i}$
for $i>l.$ Then compute $h_{1}^{-1}e_{1}\wedge ...\wedge
h_{1}^{-1}e_{d}=\det (h_{1}^{-1})\underline{e}=\mu ^{-1}\underline{v}\wedge 
\underline{w}$ where $\underline{e}=e_{1}\wedge ...\wedge e_{d}$, $%
\underline{v}=v_{1}\wedge ...\wedge v_{l}$ and $\underline{w}=e_{l+1}\wedge
...\wedge e_{d}.$ Now choose $\mu \in k$ so that $\det (h_{1})=1,$ then $%
|\mu |_{k}=||\underline{v}\wedge \underline{w}||=d(V_{a}^{\omega
},H_{a}^{\omega }).$ Then $||h_{1}^{-1}||\leq |\mu ^{-1}|_{k}$ when $k$ is
ultrametric and $||h_{1}^{-1}||\leq \sqrt{d}|\mu ^{-1}|_{k}$ when $k$ is
archimedean.

So $h_{1}ah_{1}^{-1}$ stabilizes $V$ and $H$. Now applying Lemma \ref{Conj}
on $V$ and on $H$, we can find $h_{0}\in SL_{d}(k)$, stabilizing $V$ and $H$
such that $||h_{0}h_{1}ah_{1}^{-1}h_{0|V}^{-1}||\leq C_{k,1}\Lambda _{k}(a)$
and $||h_{0}h_{1}ah_{1}^{-1}h_{0|H}^{-1}||\leq C_{k,1}\lambda _{k}^{\omega
}(a)$ and $||h_{0}^{-1}||\leq ||h_{1}ah_{1}^{-1}||^{\frac{d+1}{2}}.$ Set $%
h=h_{0}h_{1}$ we are done.%
\endproof%

\subsection{Proof of Lemma \ref{prox}}

For $l\in [1,d-1]$ set as above $V=\left\langle e_{2},...,e_{l}\right\rangle 
$ and $H=\left\langle e_{l+1},...,e_{d}\right\rangle .$ Let $b\in SL_{d}(k)$
be such that $bV=V$ and $bH=H$ and let $\pi _{b}$ be the linear projection
onto $V$ with kernel $H.$ We first claim that for every $u\in \Bbb{P}(k^{d})$%
\begin{equation}
d(u,H)\cdot d(bu,\pi _{b}(bu))\leq ||b_{|H}||\cdot ||(b_{|V})^{-1}||
\label{cont0}
\end{equation}
Indeed, note that $d(u,H)=\frac{||u\wedge \underline{w}||}{||u||}\leq \frac{%
||\pi _{b}(u)||}{||u||}$ where $\underline{w}=e_{l+1}\wedge ...\wedge e_{d}$
and writing $u=\pi _{b}(u)+\pi _{b}(u)^{\bot }$ we have $d(bu,\pi _{b}(bu))=%
\frac{||b_{|H}\pi _{b}(u)^{\bot }\wedge \pi _{b}(u)||}{||bu||}\leq \frac{%
||b_{|H}||\cdot ||u||}{||b\pi _{b}(u)||}.$ Combining both inequalities we
get $(\ref{cont0})$.

Now we claim that for any $u\neq v\in \Bbb{P}(k^{d})$ we claim that: 
\begin{equation}
\frac{d(u,H)\cdot d(v,H)\cdot d(bu,bv)}{d(u,v)}\leq \max \{||\Lambda
^{2}b_{|V}||,||\Lambda ^{2}b_{|H}||,||b_{|H}||\cdot ||b_{|V}||\}\cdot
||(b_{|V})^{-1}||^{2}  \label{lipa0}
\end{equation}
Indeed, using Cartan's $\Bbb{K}_{k}A\Bbb{K}_{k}$ decomposition on $V$ and $H$
separately, we may assume that $b$ is diagonal $diag(\alpha _{1},...,\alpha
_{l})$. Then write $bu\wedge bv=bu_{|V}\wedge bv_{|V}+bu_{|H}\wedge
bv_{|V}+bu_{|V}\wedge bv_{|H}+bu_{|H}\wedge bv_{|H}.$ Since $bu_{|H}\wedge
bv_{|V}+bu_{|V}\wedge bv_{|H}=\sum_{1\leq i\leq l}\alpha _{i}e_{i}\wedge
(u_{i}v_{|H}-v_{i}u_{|H})$ we get $||bu\wedge bv||\leq \max \{||\Lambda
^{2}b_{|V}||,||\Lambda ^{2}b_{|H}||,||b_{|H}||\cdot ||b_{|V}||\}\cdot
||u\wedge v||.$ On the other hand $||bu||\geq \frac{||u_{|V}||}{%
||b_{|V}^{-1}||}$ (resp. $||bv||\geq \frac{||v_{|V}||}{||b_{|V}^{-1}||}$)
and $d(u,H)\leq \frac{||u_{|V}||}{||u||}$ (resp.$d(v,H)\leq \frac{||v_{|V}||%
}{||v||}$). This shows (\ref{lipa0}).

We now prove (\ref{cont}) and (\ref{lipa}). For $n\in \Bbb{N}$ and $a\in
SL_{d}(k)$, we may apply Lemma \ref{interm2} and Remark \ref{Conjplane} to $%
a $ and get $h\in SL_{d}(k)$ with $\widetilde{a}=hah^{-1}\in SL_{d}(k)$
satisfying $\widetilde{a}V=V$ (resp. $\widetilde{a}H=H$) and $||\widetilde{a}%
_{|H}||\leq C_{k,1}\lambda _{k}^{\omega }(a)$ (resp. $||\widetilde{a}%
_{|V}||\leq C_{k,1}\Lambda _{k}(a)$ and $||\Lambda ^{2}\widetilde{a}%
_{|V}||\leq C_{k,1}^{4}\Lambda _{k}(a)\lambda _{k}(a)$)$.$ Note that $%
\left\| \widetilde{a}_{|V}^{-1}\right\| \leq \frac{1}{\det \widetilde{a}_{|V}%
}\left\| \widetilde{a}_{|V}\right\| ^{l_{\omega }-1}\leq \frac{1}{\Lambda
_{k}^{\omega }(a)}\cdot \left( C_{k,1}\cdot \frac{\Lambda _{k}(a)}{\Lambda
_{k}^{\omega }(a)}\right) ^{l_{\omega }-1}.$ Let $b=(\widetilde{a})^{n}.$
Then (\ref{cont0}) and (\ref{lipa0}) translate as 
\begin{equation*}
d(u,H_{a}^{\omega })\cdot d(a^{n}u,\pi _{a}^{\omega }(a^{n}u))\leq
Lip(h^{-1})^{2}\cdot \left( C_{k,1}^{l_{\omega }}\cdot \left( \frac{\Lambda
_{k}(a)}{\Lambda _{k}^{\omega }(a)}\right) ^{l_{\omega }-1}\cdot \frac{%
\lambda _{k}^{\omega }(a)}{\Lambda _{k}^{\omega }(a)}\right) ^{n}
\end{equation*}
\begin{equation*}
d(u,H_{a}^{\omega })\cdot d(v,H_{a}^{\omega })\cdot d(a^{n}u,a^{n}v)\leq
Lip(h^{-1})^{3}Lip(h)\cdot \left( C_{k,1}^{2l_{\omega }+2}\cdot \left( \frac{%
\Lambda _{k}(a)}{\Lambda _{k}^{\omega }(a)}\right) ^{2l_{\omega }-1}\cdot 
\frac{\lambda _{k}(a)}{\Lambda _{k}^{\omega }(a)}\right) ^{n}
\end{equation*}
Recall that by Lemma \ref{lip}, $Lip(h)$ and $Lip(h^{-1})$ are at most $%
||h^{-1}||^{2d}.$ From Lemma \ref{interm} we have $d(V_{a}^{\omega
},H_{a}^{\omega })^{-1}\leq \left( C_{k}\cdot L_{k}(a)\cdot ||a||^{l}\right)
^{\left( \binom{d}{l}-1\right) }$. Then from Lemma \ref{interm2} 
\begin{equation*}
||h^{-1}||\leq \left( \sqrt{C_{k,1}}\cdot (L_{k}(a)C_{k})^{\left( \binom{d}{l%
}-1\right) }\right) ^{\frac{d(d+1)}{2}}\left\| a\right\| ^{\frac{d-1}{2}%
+l\left( \binom{d}{l}-1\right) \frac{d(d+1)}{2}}
\end{equation*}
Note that, when $k$ is archimedean, inequality (\ref{cont}) is trivial if $%
\lambda _{k}^{\omega }(a)\geq \frac{1}{2}\Lambda _{k}^{\omega }(a).$
Therefore we may assume in the archimedean case that $L_{k}(a)\leq 2$. As $%
8d(\frac{d-1}{2}+l\left( \binom{d}{l}-1\right) \frac{d(d+1)}{2})\leq 10^{d}$
estimating the constant we do indeed obtain (\ref{lipa}) and (\ref{cont}). 
\endproof%

\section{Ping-Pong\label{PP}}

In this technical section, we work with a fixed local field $k$ and we
explain how to construct two short words $x$ and $y$ with letters in some
finite set $F$ in $SL_{d}(k)$ such that $x$ and $y$ form a ping-pong pair
and thus generate a free subgroup. The goal of this introductory paragraph
is to give a list of several conditions of geometric nature $(i)$ to $(vi)$
on $F$ and state two lemmas, Lemma \ref{veryprox} and \ref{pair} below,
which assert precisely that these conditions are sufficient to construct the
ping-pong pair. These two statements are the only ones which will be used in
further sections.

As in Tits \cite{Tits}, the construction of the ping-pong pair follows two
steps. First, starting from a proximal element $a$ lying in $F$ or in a
bounded power of $F,$ we need to build a short word with letters in $F$, say 
$x,$ such that both $x$ and $x^{-1}$ are proximal elements (Lemma \ref
{veryprox}). Second, we need to find a conjugate of it, say $y=cxc^{-1}$
such that $x$ and $y$ together play ping-pong (Lemma \ref{pair}).

The construction presented here follows verbatim that of Tits. But while
Tits needed only asymptotic statements which held for sufficienlty high
powers of group elements, no matter how high, we need to have control on the
length of the words. We thus have to give a quantified version of Tits'
argument and give precise estimates at each step. More importantly, while
Tits did not need to care about the choice of a distance on $\Bbb{P}(k^{d})$
(any one inside the ``admissible'' class he defined was good for his
purposes), it is crucial for us that we work with the Fubini-Study distance
introduced in Section \ref{proximality}. The reason is that all constants
then disappear and are equal to $1$ for all ultrametric local fields, hence
giving to us the possibility to bound the length of the generators of the
free subgroup \textbf{independently of the choice of the local field.}

Let $(k_{i})_{1\leq i\leq 4}$ be four positive integers and $\varepsilon
_{0},T_{0},T_{1},T_{2}>0$ be positive real numbers. Let $\varepsilon >0$
with $\varepsilon \leq \varepsilon _{0}/12d^{2}.$ Let $k_{0}$ be a local
field. Suppose $F\subset SL_{d}(k_{0})$ is a finite set containing $1.$ All
eigenvalues and eigenspaces of elements in the group generated by $F$ are
defined over a fixed finite extension of $k_{0}$ of degree at most $d!$. Let 
$k$ be this extension. For a subspace $V$ in $k^{d}$ we denote by $V^{\bot }$
its orthogonal in the dual space of $k^{d}.$ We say that a non-trivial
subspace $W$ of $k^{d}$ is $a$\textit{-admissible} for $a\in SL_{d}(k_{0})$
if it is a sum of generalized eigenspaces of $a$. We also denote by $W^{c}$
its complementary, i.e. the sum of the remaining generalized eigenspaces, so
that $k^{d}=W\oplus W^{c}$.

\textbf{List of Conditions for ping-pong (i)-(vi):}

\noindent Let $a\in F^{k_{1}},$ $b\in F^{k_{2}},$ $t\in F^{k_{3}}.$ Assume

$\mathbf{(i)}$ $a$ is proximal

$\mathbf{(ii)}$%
\begin{equation}
||F||_{k}>C_{k,1}^{2d}  \label{a1b}
\end{equation}

$\mathbf{(iii)}$%
\begin{equation}
\left( \frac{\Lambda _{k}(a)}{\lambda _{k}(a)}\right) ^{\frac{1}{\varepsilon
_{0}}}\geq \Lambda _{k}(a)\geq ||F||_{k}^{T_{1}}  \label{a1c}
\end{equation}

$\mathbf{(iv)}$ For any $a$-admissible subspace $W$ (see Def. \ref{admiss})
we have 
\begin{equation}
d(^{t}b^{\pm 1}\cdot H_{a}^{\bot },W^{\bot })>||F||_{k}^{-T_{0}}  \label{bb}
\end{equation}

$\mathbf{(v)}$ For any $a$-admissible subspace $W$ we have 
\begin{eqnarray}
d(tV_{a},W^{c}+W\cap b^{-1}H_{a}) &\geq &||F||_{k}^{-T_{0}}  \label{u1} \\
d(t^{-1}V_{a},W^{c}+W\cap bH_{a}) &\geq &||F||_{k}^{-T_{0}}
\end{eqnarray}

Note that condition $(\ref{bb})$ on $b$ implies that $W^{c}+W_{1}\cap b^{\pm
1}H_{a}$ are hyperplanes, so these distances are computable via $(\ref
{disttoH})$.

\begin{lemma}
\label{veryprox}There is $\tau _{1}(d,\varepsilon )\in \Bbb{N}$ and $\tau
_{3}=\tau _{3}(d,k_{1},k_{2},k_{3},\varepsilon _{0},\varepsilon
,T_{0},T_{1})\in \Bbb{N}$ such that if $T_{1}\geq \tau _{1}$ and $T_{3}\geq
\tau _{3}$, there is $l=l(d,k_{1},k_{2},k_{3},\varepsilon _{0},\varepsilon
,T_{0},T_{1},T_{3})\in \Bbb{N}$ such that for some $l_{0},l_{1}\in [0,l]$
the element $x=a^{l_{0}}ba^{-l_{1}}t$ is proximal as well as $x^{-1}$ and $%
\Lambda _{k}(x)\geq \Lambda _{k}(a)^{T_{3}}\lambda _{k}(x)$ and $\Lambda
_{k}(x^{-1})\geq \Lambda _{k}(a)^{T_{3}}\lambda _{k}(x^{-1})$ and $%
d(V_{x},H_{x})\geq \Lambda _{k}(a)^{-2T_{3}}$ and $d(V_{x^{-1}},H_{x^{-1}})%
\geq \Lambda _{k}(a)^{-2T_{3}}.$
\end{lemma}

We let $k_{4}=2k_{1}l+k_{2}+k_{3}$ so that $x\in F^{k_{4}}.$

\begin{remark}
\label{benoist}As Y. Benoist observed in \cite{Ben} (see also J-F. Quint 
\cite{Quin}) it is possible to construct Zariski-dense semi-groups, say in $%
SL_{3}(\Bbb{Q}_{p})$ which are made exclusively of proximal elements whose
inverses are not proximal. Hence our method does not allow in general (the $%
SL_{2}$ case is fine however) to construct the generators of a free subgroup
as positive words in $F$.
\end{remark}

Assume $T_{1}$ and $T_{3}$ satisfy the assumptions of Lemma \ref{veryprox}
let $x$ be the element we get. Assume that there is $c\in F^{k_{5}}$ such
that

$\mathbf{(vi)}$%
\begin{equation}
d(c^{\pm 1}(V_{x}\cup V_{x^{-1}}),H_{x}\cup H_{x^{-1}})\geq \frac{1}{%
||F||_{k}^{T_{2}}}  \label{c1}
\end{equation}

\begin{lemma}
\label{pair}Then there is $l_{2}=l_{2}(d,(k_{i})_{1\leq i\leq 5},\varepsilon
,\varepsilon _{0},(T_{i})_{0\leq i\leq 3})\in \Bbb{N}$ such that for every $%
n\geq l_{2},$ $x^{n}$ and $y=cx^{n}c^{-1}$ play ping-pong on $\Bbb{P}(k^{d})$
and generate a free subgroup of $SL_{d}(k).$
\end{lemma}

\subsection{Cayley-Hamilton trick}

In \cite{Tits} Tits used the fact that if $a_{0}\in GL_{d}(k)$ has all
eigenvalues of the same modulus and if a vector $v$ lies far from a
hyperplane $H$ then for a set of positive density of $n\in \Bbb{N}$ the
vectors $a^{n}\cdot v$ lie far from $H$. In \cite{EMO}, Eskin-Mozes-Oh made
clever use of the Cayley-Hamilton theorem in order to show a statement of a
similar nature which also gives a bound on the smallest appropriate $n$. The
following lemma is a reformulation of the same trick.

Recall that $C_{k,2}$ is $d2^{d}$ when $k$ is archimedean and $1$ when $k$
is ultrametric. The following lemma, which we will use in the proof of Claim
0 below, expresses the same idea.

\begin{lemma}
\label{Cayley}Let $a_{0}\in GL_{d}(k)$, let $H$ be a hyperplane in $k^{d}$
and let $v\in \Bbb{P}(k^{d}).$ Then there is some integer $j_{0}\in [1,d-1]$
such that 
\begin{equation}
d(a_{0}^{j_{0}}v,H)\geq \frac{1}{C_{k,2}}\cdot \left( \frac{\Lambda
_{k}(a_{0})}{||a_{0}||_{k}}\right) ^{j_{0}}\cdot \frac{|\det a_{0}|_{k}}{%
\Lambda _{k}(a_{0})^{d}}\cdot d(v,H)  \label{simil}
\end{equation}
\end{lemma}

\proof%
Let $\Lambda \in k$ such that $|\Lambda |=\Lambda _{k}(a_{0})$ and set $%
\widetilde{a_{0}}=\frac{a_{0}}{\Lambda }.$ According to the Cayley-Hamilton
theorem, there are coefficients $(c_{j})_{1\leq j\leq d-1}$ in $k$ such that 
$\sum_{j=1}^{d-1}c_{j}\widetilde{a_{0}}^{j}=\det \widetilde{a_{0}}.$
Moreover $|c_{j}|_{k}\leq \binom{d}{j}\leq 2^{d}$ when $k$ is archimedean,
when $|c_{j}|_{k}\leq 1$ when $k$ is ultrametric. Let $f$ be a linear form
on $k^{d}$ with $||f||_{k}=1$ and $\ker f=H.$ There must exist some $%
j_{0}\in [1,d-1]$ such that $|c_{j_{0}}f(\widetilde{a_{0}}%
^{j_{0}}v)|_{k}\geq \frac{|\det \widetilde{a_{0}}|_{k}}{C_{k,1}}\cdot
|f(v)|_{k}$ where $C_{k,1}$ is $d$ if $k$ is archimedean and $1$ if $k$ is
ultrametric. Hence $|f(a_{0}^{j_{0}}v)|_{k}\geq \frac{1}{C_{k,2}}\cdot
\Lambda _{k}(a_{0})^{j_{0}}\cdot \frac{|\det a_{0}|_{k}}{\Lambda
_{k}(a_{0})^{d}}\cdot |f(v)|_{k}$ and $(\ref{simil})$ follows. 
\endproof%

\subsection{Proof of Lemma \ref{veryprox}}

Recall that $a$ is proximal but $a^{-1}$ may not be. However, as we have
fixed $\varepsilon >0$, Lemma \ref{almostprox} gives us some $\omega $ for
which $a^{-1}$ is almost proximal. Let $\alpha =\frac{\lambda _{k}^{\omega
}(a^{-1})}{\Lambda _{k}^{\omega }(a^{-1})}.$ It also give $\eta =\eta
(d,\varepsilon )>0.$ Recall that $\varepsilon _{0}$, $T_{0}$ and $T_{1}$ are
defined in $(\ref{a1b})$ to $(\ref{u1})$. We assume here that $T_{1}\geq
\tau _{1}:=\max \{2/\eta ,3/\eta \varepsilon ,4/\varepsilon _{0}\}$ and $%
\varepsilon \leq \varepsilon _{0}/12d^{2}$ and let $\varepsilon _{1}=\frac{%
\varepsilon _{0}}{4}$.

\textbf{Claim 0:} There is $n_{0}=n_{0}(k_{1},T_{0},T_{1},\varepsilon ,d)\in 
\Bbb{N}$ such that for all $n\geq n_{0}$ there exists $j_{0}(n),j_{1}(n)\in
[1,d-1]$ such that 
\begin{equation*}
\min \left\{
d(a^{-nj_{0}(n)}t^{-1}V_{a},bH_{a}),d(a^{-nj_{1}(n)}tV_{a},b^{-1}H_{a})%
\right\} \geq \frac{||F||_{k}^{-r}}{C_{k,2}\cdot (C_{k,1}\cdot \alpha
^{-\varepsilon })^{dn}\cdot C_{k}^{p(d)+1}}
\end{equation*}
where $r=2T_{0}+k_{1}(2p(d)+d(d-1))$ and $p(d)=10^{d}.$

\textbf{Claim 1:} Let $\varepsilon _{1}=\frac{\varepsilon _{0}}{4}$. There
exists $n_{1}=n_{1}(k_{1},k_{2},k_{3},T_{0},T_{1},\varepsilon ,\varepsilon
_{0},d)\in \Bbb{N}$ such that for every $n\geq n_{1}$ the ball $%
B_{n}=B(V_{a},\Lambda _{k}(a)^{-\varepsilon _{1}n})$ (resp. $B_{n}^{\prime
}=B(t^{-1}V_{a},\Lambda _{k}(a)^{-\varepsilon _{1}n}$) is mapped into $%
B_{n}^{-}=B(V_{a},\Lambda _{k}(a)^{-2\varepsilon _{1}n})$ (resp. $%
B_{n}^{\prime -}=B(t^{-1}V_{a},\Lambda _{k}(a)^{-2\varepsilon _{1}n})$) by $%
x_{n}=a^{nj_{0}(n)}ba^{-nj_{1}(n)}t$ (resp. $x_{n}^{-1}$).

\textbf{Claim 2:} Under the assumptions of Claim 1, there is $n_{4}\in \Bbb{N%
}$ depending only on $k_{1},k_{2},k_{3},T_{0},T_{1},\varepsilon ,\varepsilon
_{0}$ and $d$ such that for any $n\geq n_{4}$ we also have $%
Lip(x_{n|B_{n}})\leq \Lambda _{k}(a)^{-\varepsilon _{1}n}$ (resp. $%
Lip(x_{n|B_{n}^{-}}^{-1})\leq \Lambda _{k}(a)^{-\varepsilon _{1}n}$).

The proofs of these claims are straightforward once we have at our disposal
the Lemmas proved in Section \ref{proximality} and in particular Lemma \ref
{prox}. Nevertheless we provide full details in the next paragraph below.

With these claims in hands we can quickly prove Lemma \ref{veryprox}. Indeed
let $n=T_{3}/\varepsilon _{1}.$ If $T_{3}\geq \varepsilon _{1}\cdot \max
\{n_{0},n_{1},n_{4}\}$ we get by Claim 2 and 3 that $x_{n}$ sends $B_{n}$
into itself and $x_{n}^{-1}$ sends $B_{n}^{-}$ into itself, while the
Lipschitz constants are $\leq \Lambda _{k}(a)^{-\varepsilon _{1}n}.$ We are
thus in a position to apply Tits Converse Lemma, Lemma \ref{TitsConv}, which
says that $x_{n}$ and $x_{n}^{-1}$ are proximal and satisfy $\frac{\lambda
_{k}(x_{n})}{\Lambda _{k}(x_{n})},\frac{\lambda _{k}(x_{n}^{-1})}{\Lambda
_{k}(x_{n}^{-1})}\leq \Lambda _{k}(a)^{-T_{3}}$. Finally by Claim 2, $x_{n}$
maps $B_{n}$ into the smaller ball $B_{n}^{-}$, which must then contain $%
V_{x_{n}}$ while $B_{n}$ cannot intersect $H_{x_{n}}.$ It follows that $%
d(V_{x_{n}},H_{x_{n}})\geq d(B_{n},(B_{n}^{-})^{c})$. But we see that in
both the archimedean and the ultrametric case: 
\begin{equation*}
d(B_{n},(B_{n}^{-})^{c})\geq \frac{1}{C_{k}}\Lambda _{k}(a)^{-T_{3}}\geq
\Lambda _{k}(a)^{-2T_{3}}
\end{equation*}
as soon as $\Lambda _{k}(a)^{-T_{3}}<1/C_{k},$ which holds if $T_{3}\geq 1$
for instance. A similar argument takes place for $x_{n}^{-1}.$ This ends the
proof of Lemma \ref{veryprox}.

\subsubsection{Proof of Claim 0}

Let $W=V_{a^{-1}}^{\omega }$ and hence $W^{c}=H_{a^{-1}}^{\omega }$ and $\pi 
$ the projection on $W$ with kernel $W^{c}.$ Recall that $(\ref{a1c})$ gives 
\begin{equation}
\alpha \leq \Lambda _{k}(a)^{-\eta }\leq ||F||^{-\eta T_{1}}  \label{upalph}
\end{equation}
When $k$ is archimedean this together with $(\ref{a1b})$ forces $%
L_{k}^{\omega }(a^{-1})\leq 2$ if $\eta T_{1}\geq 2$. Indeed, we have $%
\alpha \leq ||F||^{-\eta T_{1}}\leq \frac{1}{2}$ i.e. $\lambda _{k}^{\omega
}(a^{-1})\leq \frac{1}{2}\Lambda _{k}^{\omega }(a^{-1}).$

We do the calculation for $u=tV_{a}$, keeping in mind that an entirely
analogous calculation can be done for $u^{-}=t^{-1}V_{a}$ at the same time
at each step. Let $n\in \Bbb{N}$ be arbitrary.

Since as $d(u,W^{c})\geq d(u_{0},W^{c}+W\cap b^{-1}H_{a})\geq ||F||^{-T_{0}}$
we can combine Lemmas \ref{almostprox}, \ref{prox} $(a)$ and $(\ref{u1})$ to
get 
\begin{eqnarray}
d(a^{-n}u,\pi (a^{-n}u)) &\leq &(C_{k}\cdot ||a^{-1}||_{k})^{p(d)}\cdot
\left( C_{k,1}^{l_{\omega _{1}}}\cdot \left( \frac{\Lambda _{k}(a^{-1})}{%
\Lambda _{k}^{\omega }(a^{-1})}\right) ^{l_{\omega }-1}\cdot \frac{\lambda
_{k}^{\omega }(a^{-1})}{\Lambda _{k}^{\omega }(a^{-1})}\right) ^{n}\cdot
d(u,W^{c})^{-1}  \label{other0} \\
&\leq &C_{k}^{p(d)}\cdot ||F||_{k}^{T_{0}+dk_{1}p(d)}\cdot \left(
C_{k,1}^{d}\cdot \alpha ^{1-d\varepsilon }\right) ^{n}  \notag
\end{eqnarray}
On the other hand according to Lemma \ref{dist1}, 
\begin{equation}
d(\pi u,W\cap b^{-1}H_{a})\geq d(u,W^{c}+W\cap b^{-1}H_{a})\cdot d(W,W^{c})
\label{other1}
\end{equation}
But by Lemma \ref{interm} 
\begin{equation}
d(W,W^{c})^{-1}\leq (C_{k}^{2}\cdot ||a^{-1}||^{l_{\omega }})^{\left( \binom{%
d}{l_{\omega }}-1\right) }\leq (C_{k}\cdot ||a^{-1}||_{k})^{p(d)}
\label{other2}
\end{equation}
because when $k$ is archimedean $L_{k}^{\omega }(a^{-1})\leq 2$ as explained
above$.$ Hence $(\ref{other1})$ $(\ref{other2})$ and $(\ref{u1})$ give 
\begin{equation}
d(\pi u,W\cap b^{-1}H_{a})^{-1}\leq C_{k}^{p(d)}\cdot
||F||_{k}^{T_{0}+dk_{1}p(d)}  \label{other3}
\end{equation}
We may now apply Lemma \ref{Cayley} to $a_{0}=a^{-n}$ restricted to $W.$ We
find $j_{1}\in [1,d-1]$ such that 
\begin{equation}
d(a^{-nj_{1}}\pi u,W\cap b^{-1}H_{a})^{-1}\leq C_{k,2}\cdot \left( \frac{%
||a^{-1}||_{k}}{\Lambda _{k}(a^{-1})}\right) ^{nj_{1}}\cdot \left( \frac{%
\Lambda _{k}(a^{-1})}{\Lambda _{k}^{\omega }(a^{-1})}\right) ^{\ell _{\omega
}n}\cdot d(\pi u,W\cap b^{-1}H_{a})^{-1}  \label{other4}
\end{equation}
But Lemma \ref{Conj} applied to $a^{-1}$ gives an $h\in SL_{d}(k)$ such that 
$||ha^{-1}h^{-1}||\leq C_{k,1}\cdot \Lambda _{k}(a^{-1})$ and $\max
\{||h||,||h^{-1}||\}\leq ||a^{-1}||^{\frac{d-1}{2}}.$ Hence $%
||a^{-n}||_{k}\leq ||h||\cdot ||h^{-1}||\cdot (C_{k,1}\cdot \Lambda
_{k}(a^{-1}))^{n}$ and $\frac{||a^{-n}||_{k}}{\Lambda _{k}(a^{-n})}\leq
||a^{-1}||^{d-1}\cdot C_{k,1}^{n}.$ Thus combining $(\ref{other3})$ and $(%
\ref{other4})$ and bearing in mind that 
\begin{equation}
\frac{\Lambda _{k}(a^{-1})}{\Lambda _{k}^{\omega }(a^{-1})}\leq \alpha
^{-\varepsilon }  \label{almostproxupbnd}
\end{equation}
(this is \ref{almostproxineq}), we get 
\begin{equation}
d(a^{-nj_{1}}\pi u,W\cap b^{-1}H_{a})^{-1}\leq C_{k,2}\cdot (C_{k,1}\alpha
^{-\varepsilon })^{nd}\cdot C_{k}^{p(d)}\cdot
||F||_{k}^{T_{0}+dk_{1}p(d)+k_{1}d^{2}(d-1)}
\end{equation}
Compare (\ref{other0}) and (\ref{other5}). When $k$ is ultrametric 
\begin{eqnarray*}
d(a^{-n}u,\pi (a^{-n}u)) &\leq &||F||_{k}^{T_{0}+dk_{1}p(d)}\alpha
^{n(1-d\varepsilon )} \\
&<&||F||_{k}^{-T_{0}-dk_{1}p(d)-k_{1}d^{2}(d-1)}\alpha ^{\varepsilon nd} \\
&<&d(a^{-nj_{1}}\pi u,W\cap b^{-1}H_{a})
\end{eqnarray*}
as soon as $\alpha ^{n(1-2d\varepsilon )}<||F||_{k}^{-r}$ where $%
r=r(T_{0},d):=2T_{0}+dk_{1}(2p(d)+d(d-1)).$ As $\alpha ^{-1}\geq \Lambda
_{k}(a)^{\eta }$ by Lemma \ref{almostprox}, this happens as soon as 
\begin{equation*}
n>\frac{r}{T_{1}\eta (1-\varepsilon d)}
\end{equation*}
Similarly, if $k$ is archimedean, $d(a^{-n}u,\pi (a^{-n}u))\leq \frac{1}{2}%
d(a^{-nj_{1}}\pi u,W\cap b^{-1}H_{a})$ as soon as $n>n_{0}(T_{0},T_{1},%
\varepsilon ,d)$ for some computable constant $n_{0}.$ Finally whether $k$
is archimedean or ultrametric we get: 
\begin{equation}
d(a^{-nj_{1}}u,W\cap b^{-1}H_{a})\geq C_{k,2}^{-1}\cdot (C_{k,1}^{-1}\cdot
\alpha ^{\varepsilon })^{dn}\cdot C_{k}^{-p(d)-1}\cdot
||F||_{k}^{-T_{0}-k_{1}p(d)-k_{1}d^{2}(d-1)}\cdot  \label{other5}
\end{equation}
Finally applying Lemma \ref{dist0} and $(\ref{bb})$ we get 
\begin{eqnarray*}
d(a^{-nj_{1}}u,b^{-1}H_{a}) &\geq &d(a^{-nj_{1}}u,W\cap b^{-1}H_{a})\cdot
d(^{t}b^{-1}\cdot H_{a}^{\bot },W^{\bot }) \\
&\geq &C_{k,2}^{-1}\cdot (C_{k,1}^{-1}\cdot \alpha ^{\varepsilon
})^{dn}\cdot C_{k}^{-p(d)-1}\cdot ||F||_{k}^{-r}
\end{eqnarray*}
This ends the proof of Claim 0.

\subsubsection{Proof of Claim 1}

First recall as in Claim 0 that $L_{k}^{\omega }(a^{-1})\leq 2$ when $k$ is
archimedean (since $\eta T_{1}\geq 2,$ which we assume). We give the proof
for $x_{n}$ and $B_{n}$ keeping in mind that the same arguments are being
performed at the same time and at each step for $x_{n}^{-1}$ and $B_{n}^{-}.$

We first justify the following:

\textbf{Claim 1.1:} There is $m_{0}=m_{0}(d,\varepsilon
,k_{3},T_{0},T_{1})\in \Bbb{N}$ such that for $n\geq m_{0}$ and $u\in
B(V_{a},\alpha ^{3d\varepsilon n})$%
\begin{equation}
d(tu,W^{c})\geq \frac{1}{C_{k}}d(tV_{a},W^{c})\geq \frac{1}{%
C_{k}||F||_{k}^{T_{0}}}  \label{tuaway}
\end{equation}
Indeed, the second inequality is just $(\ref{u1}),$ while to get the first,
it is enough that $d(tu,tV_{a})<\frac{1}{C_{k}||F||_{k}^{T_{0}}}\leq \frac{1%
}{C_{k}}d(tV_{a},W^{c})$ (recall that $C_{k}$ is $2$ is $k$ is archimedean
and $1$ if $k$ is ultrametric). But 
\begin{equation*}
d(tu,tV_{a})\leq Lip(t)\cdot d(u,V_{a})\leq ||F||^{2dk_{3}}\alpha
^{3d\varepsilon n}
\end{equation*}
Thus the existence of $m_{0}$ follows from $(\ref{upalph})$ $(\ref{a1b})$
and $(\ref{a1c}).$ Hence $(\ref{tuaway}).$

\textbf{Claim 1.2.:} For some $n_{2}=n_{2}(\varepsilon
,d,T_{0},T_{1},k_{1},k_{3})\in \Bbb{N}$ and all $n\geq n_{2}$ we have 
\begin{eqnarray}
d(a^{-nj_{1}(n)}tu,b^{-1}H_{a}) &\geq &\frac{1}{C_{k}}%
d(a^{-nj_{1}(n)}tV_{a},b^{-1}H_{a})  \label{other5bis} \\
&\geq &C_{k,2}^{-1}\cdot (C_{k,1}^{-1}\cdot \alpha ^{\varepsilon
})^{dn}\cdot C_{k}^{-p(d)-2}\cdot ||F||_{k}^{-r}  \notag
\end{eqnarray}
for all $u\in B(V_{a},\alpha ^{3d\varepsilon n}).$

\textit{Proof of Claim 1.2.:} Indeed, to show this it is enough that 
\begin{equation}
d(a^{-nj_{1}(n)}tu,a^{-nj_{1}(n)}tV_{a})<\frac{1}{C_{k}}%
d(a^{-nj_{1}(n)}tu,b^{-1}H_{a}),  \label{other6}
\end{equation}
which by Claim 0 reduces to show 
\begin{equation}
d(a^{-nj_{1}(n)}tu,a^{-nj_{1}(n)}tV_{a})<C_{k,2}^{-1}\cdot
(C_{k,1}^{-1}\cdot \alpha ^{\varepsilon })^{dn}\cdot C_{k}^{-p(d)-2}\cdot
||F||_{k}^{-r}  \label{other7}
\end{equation}
But bearing in mind $(\ref{almostproxupbnd})$ Lemma \ref{prox} $(\ref{lipa}%
), $ we have for $n\geq m_{0}$

\begin{eqnarray*}
\frac{d(a^{-nj_{1}(n)}tu,a^{-nj_{1}(n)}tV_{a})}{d(u,V_{a})} &\leq
&Lip(t)\cdot \frac{d(a^{-nj_{1}(n)}tu,a^{-nj_{1}(n)}tV_{a})}{d(tu,tV_{a})} \\
&\leq &Lip(t)\cdot (C_{k}^{2}\cdot ||a^{-1}||)^{p}\cdot \frac{\left(
C_{k,1}^{2d+2}\alpha ^{-\varepsilon (2d-1)}\right) ^{n}}{d(tu,W^{c})\cdot
d(tV_{a},W^{c})} \\
&\leq &C_{k}^{2p+1}\cdot ||F||_{k}^{k_{1}dp+k_{3}2d+2T_{0}}\cdot \left(
C_{k,1}^{2d+2}\alpha ^{-\varepsilon (2d-1)}\right) ^{n}
\end{eqnarray*}
Hence we get $(\ref{other7})$ as soon as 
\begin{equation*}
C_{k}^{3p+3}\cdot C_{k,2}\cdot ||F||_{k}^{k_{1}dp+k_{3}2d+2T_{0}+r}\cdot
\left( C_{k,1}^{3d+2}\alpha ^{-\varepsilon (3d-1)}\right) ^{n}\cdot
d(u,V_{a})<1
\end{equation*}
Since $u\in B(V_{a},\alpha ^{3d\varepsilon n})$ this holds as soon as 
\begin{equation}
C_{k}^{3p+3}\cdot C_{k,2}\cdot ||F||_{k}^{k_{1}dp+k_{3}2d+2T_{0}+r}\cdot
\left( C_{k,1}^{3d+2}\alpha ^{\varepsilon }\right) ^{n}<1  \label{other8}
\end{equation}
Since $\alpha ^{\varepsilon }\leq ||F||_{k}^{-\eta \varepsilon T_{1}}$ by $(%
\ref{upalph})$ and $C_{k,1}^{3d+2}\leq ||F||_{k}^{2}$ by $(\ref{a1b})$ while
we assumed $\eta \varepsilon T_{1}\geq 3,$ we get the existence of $%
n_{2}=n_{2}(\varepsilon ,d,T_{0},T_{1},k_{1},k_{3})\in \Bbb{N}$ for which (%
\ref{other8}) holds for $n\geq n_{1}.$ Hence (\ref{other5bis}) holds and
Claim 1.2. is proved.

With (\ref{other5bis}) in hand we can apply Lemma \ref{prox} $(\ref{cont})$
to positive powers of $a$ this time and get:

\textbf{Claim 1.3.:} Suppose $\varepsilon _{0}\geq 12\varepsilon d^{2}$ and
fix $\varepsilon _{1}=\varepsilon _{0}/4.$ There is $n_{3}\in \Bbb{N}$
depending on $\varepsilon ,\varepsilon _{0},d,T_{0},T_{1},k_{1},k_{2},k_{3}$
such that for $n\geq n_{3}$ and $u\in B(V_{a},\Lambda _{k}(a)^{-\varepsilon
_{1}n})$ we have for $x_{n}=a^{nj_{0}(n)}ba^{-nj_{1}(n)}t$ 
\begin{equation*}
d(x_{n}u,V_{a})<\Lambda _{k}(a)^{-\varepsilon _{1}n}
\end{equation*}

\textit{Proof of Claim 1.3.: }First note that $\Lambda _{k}(a)^{-\varepsilon
_{1}}\leq \alpha ^{3d\varepsilon }$ because $\alpha ^{-1}\leq \Lambda
_{k}(a)\Lambda _{k}(a^{-1})\leq \Lambda _{k}(a)^{d}$ and $\varepsilon
_{1}=\varepsilon _{0}/4\geq 3\varepsilon d^{2}.$ Lemma \ref{prox} $(\ref
{cont})$ translates as 
\begin{eqnarray*}
d(a^{nj_{0}(n)}ba^{-nj_{1}(n)}tu,V_{a}) &\leq &(C_{k}\cdot
||a||_{k})^{p(d)}\cdot \left( C_{k,1}^{4}\cdot \frac{\lambda _{k}(a)}{%
\Lambda _{k}(a)}\right) ^{n}\cdot d(ba^{-nj_{1}(n)}tu,H_{a})^{-1} \\
&\leq &C_{k}^{p}\cdot ||F||_{k}^{k_{1}p(d)+2dk_{2}}\cdot \left( \frac{%
C_{k,1}^{4}}{\Lambda _{k}(a)^{\varepsilon _{0}}}\right) ^{n}\cdot
d(a^{-nj_{1}(n)}tu,b^{-1}H_{a})^{-1} \\
&\leq &C_{k}^{2p+2}C_{k,2}\cdot ||F||_{k}^{k_{1}p(d)+2dk_{2}+r}\cdot \left( 
\frac{C_{k,1}^{d+4}\cdot \alpha ^{-\varepsilon d}}{\Lambda
_{k}(a)^{\varepsilon _{0}}}\right) ^{n}
\end{eqnarray*}
where we have used successively $(\ref{a1c})$ and Claim 1.2. Now 
\begin{equation*}
\frac{C_{k,1}^{d+4}\cdot \alpha ^{-\varepsilon d}}{\Lambda
_{k}(a)^{\varepsilon _{0}}}\cdot \Lambda _{k}(a)^{\varepsilon _{1}}\leq 
\frac{C_{k,1}^{d+4}}{\Lambda _{k}(a)^{\varepsilon _{0}-\varepsilon
_{1}-\varepsilon d^{2}}}\leq \frac{C_{k,1}^{d+4}}{\Lambda
_{k}(a)^{\varepsilon _{0}/2}}
\end{equation*}
because $\alpha ^{-1}\leq \Lambda _{k}(a)\Lambda _{k}(a^{-1})\leq \Lambda
_{k}(a)^{d}$ and we have assumed $\varepsilon _{0}\geq 4\varepsilon d^{2}$
and $\varepsilon _{0}=4\varepsilon _{1}.$ Then the existence of $n_{3}$
follows from $(\ref{a1b})$ and $(\ref{a1c}).$ Thus Claim 1.3. is proved.

Working out the same three claims for $x_{n}^{-1}$ and $B_{n}^{\prime }$ in
place of $x_{n}$ and $B_{n}$ we get Claim 1.

\subsubsection{Proof of Claim 2}

We apply Lemma \ref{prox} $(\ref{lipa})$ to $a^{nj_{0}(n)}$ and points $%
ba^{-nj_{1}(n)}tu$ and $ba^{-nj_{1}(n)}tv$ for $u,v\in B_{n}.$ Recall that $%
L_{k}(a)\leq 2$ when $k$ is archimedean as $\frac{\lambda _{k}(a)}{\Lambda
_{k}(a)}\leq ||F||^{-T_{1}\varepsilon _{0}}\leq \frac{1}{2}$ by $(\ref{a1b})$
and $(\ref{a1c})$ and since $T_{1}\varepsilon _{0}\geq 1$. We get 
\begin{equation*}
\frac{d(x_{n}u,x_{n}v)}{d(u,v)}\leq (C_{k}^{2}||a||_{k})^{p(d)}\cdot \left( 
\frac{C_{k,1}^{4}\lambda _{k}(a)}{\Lambda _{k}(a)}\right) ^{n}\cdot
d(ba^{-nj_{1}(n)}tu,H_{a})^{-1}\cdot d(ba^{-nj_{1}(n)}tv,H_{a})^{-1}
\end{equation*}
Since $\Lambda _{k}(a)^{-\varepsilon _{1}}\leq \alpha ^{3d\varepsilon }$
Claim 1.2. applies and we get 
\begin{eqnarray*}
\frac{d(x_{n}u,x_{n}v)}{d(u,v)} &\leq &Lip(b^{-1})^{2}\cdot
||a||_{k}{}^{p(d)}\cdot \left( \frac{C_{k,1}^{4}\lambda _{k}(a)}{\Lambda
_{k}(a)}\right) ^{n}\cdot C_{k,2}^{2}\cdot (C_{k,1}\cdot \alpha
^{-\varepsilon })^{2dn}\cdot C_{k}^{4p(d)+4}\cdot ||F||_{k}^{2r} \\
&\leq &C_{k}^{4p+4}C_{k,2}^{2}\cdot ||F||_{k}^{4dk_{2}+2r+k_{1}p(d)}\cdot
\left( C_{k,1}^{4+2d}\frac{\alpha ^{-2d\varepsilon }}{\Lambda
_{k}(a)^{\varepsilon _{0}}}\right) ^{n}
\end{eqnarray*}
But 
\begin{equation*}
\frac{\alpha ^{-2d\varepsilon }}{\Lambda _{k}(a)^{\varepsilon _{0}}}\leq 
\frac{1}{\Lambda _{k}(a)^{2\varepsilon _{1}}}
\end{equation*}
Hence for some computable $n_{4},$ for all $n\geq n_{4}$ and $u,v\in B_{n}$%
\begin{equation*}
\frac{d(x_{n}u,x_{n}v)}{d(u,v)}\leq \Lambda _{k}(a)^{-\varepsilon _{1}n}
\end{equation*}
A similar argument proves the claim about $x_{n}^{-1}$ and $B_{n}^{\prime }.$
Thus Claim 2 is proved.

\subsection{Proof of Lemma \ref{pair}}

Let $n,k\in \Bbb{N},$ $k=T_{2}+3dk_{5}.$ Let $B_{k}(x)=B(V_{x},\Lambda
_{k}(a)^{-kT_{3}})$ and $B_{k}^{-}(x)=B(V_{x^{-1}},\Lambda
_{k}(a)^{-kT_{3}}).$ Similarly, let $B_{k}(c)=B(cV_{x},\Lambda
_{k}(a)^{-kT_{3}})$ and $B_{k}^{-}(c)=B(cV_{x^{-1}},\Lambda
_{k}(a)^{-kT_{3}}).$ Note that $d(u,c^{\pm 1}V_{x})<\frac{1}{C_{k}}d(c^{\pm
1}V_{x},H_{x})$ implies $d(u,H_{x})\geq \frac{1}{C_{k}}d(c^{\pm
1}V_{x},H_{x}).$ Hence if $u\in B_{k}(c),$ then as $k\geq T_{2},$ $\Lambda
_{k}(a)^{kT_{3}}>C_{k}||F||_{k}^{T_{2}}$ and 
\begin{equation*}
d(u,cV_{x})\leq \Lambda _{k}(a)^{-kT_{3}}<\frac{1}{C_{k}}d(cV_{x},H_{x})
\end{equation*}
Hence 
\begin{equation*}
d(u,H_{x})\geq \frac{1}{C_{k}||F||_{k}^{T_{2}}}\geq \frac{1}{%
||F||_{k}^{T_{2}+1}}
\end{equation*}
Similarly if $u\in B_{k}(c)$ then $d(u,H_{x^{-1}})\geq ||F||_{k}^{-T_{2}-1}$
and if $u\in B_{k}^{-}(c),$ then $d(u,H_{x}\cup H_{x^{-1}})\geq
||F||_{k}^{-T_{2}-1}.$ Finally: 
\begin{equation}
d(B_{k}^{-}(c)\cup B_{k}(c),H_{x}\cup H_{x^{-1}})\geq ||F||_{k}^{-T_{2}-1}
\label{other18}
\end{equation}
Similarly we check that 
\begin{equation*}
d(B_{k}^{-}(x)\cup B_{k}(x),cH_{x}\cup cH_{x^{-1}})\geq
||F||_{k}^{-T_{2}-2dk_{5}-1}
\end{equation*}
We know that for each $n\geq 1,$ $x^{n}$ maps $B_{k}$ into itself and $%
x^{-n} $ maps $B_{k}^{-}$ into itself. Similarly we check that $cx^{n}c^{-1}$
maps $B_{k}(c)$ into itself and $cx^{-n}c^{-1}$ maps $B_{k}^{-}(c)$ into
itself.

We now check that $x^{n}$ maps $B_{k}^{-}(c)\cup B_{k}(c)$ into $B_{k}$ and $%
x^{-n}$ maps $B_{k}^{-}(c)\cup B_{k}(c)$ into $B_{k}^{-}$. From Lemma \ref
{veryprox}, we have $\frac{\lambda _{k}(x)}{\Lambda _{k}(x)}\leq \Lambda
_{k}(a)^{-T_{3}}.$ By Lemma \ref{prox} $(\ref{cont})$ applied to $x,$ for $%
u\in \Bbb{P}(k^{d}),$%
\begin{eqnarray*}
d(x^{n}u,V_{x})\cdot d(u,H_{x}) &\leq &(C_{k}\cdot ||x||_{k})^{p(d)}\cdot
\left( C_{k,1}\cdot \frac{\lambda _{k}(x)}{\Lambda _{k}(x)}\right) ^{n} \\
&\leq &||F||_{k}^{p(d)(1+ld+k_{2}+k_{3})}\cdot \Lambda _{k}(a)^{-\frac{T_{3}n%
}{2}}
\end{eqnarray*}
Hence if $u\in B_{k}(c)\cup B_{k}^{-}(c),$ then $d(u,H_{x})\geq
||F||_{k}^{-T_{2}-1}$ by $(\ref{other18})$ and 
\begin{eqnarray*}
d(x^{n}u,V_{x}) &\leq &||F||_{k}^{p(d)(1+ld+k_{2}+k_{3})+T_{2}+1}\cdot
\Lambda _{k}(a)^{-\frac{T_{3}n}{2}} \\
&\leq &\Lambda _{k}(a)^{-kT_{3}}
\end{eqnarray*}
as soon as $n\geq n_{5}=n_{5}(l,d,(k_{i})_{i},(T_{i})_{i},k).$ Similarly 
\begin{eqnarray*}
d(x^{-n}u,V_{x^{-1}})\cdot d(u,H_{x^{-1}}) &\leq &(C_{k}\cdot
||x^{-1}||_{k})^{p(d)}\cdot \left( C_{k,1}\cdot \frac{\lambda _{k}(x^{-1})}{%
\Lambda _{k}(x^{-1})}\right) ^{n} \\
&\leq &||F||_{k}^{p(d)d(1+ld+k_{2}+k_{3})}\cdot \Lambda _{k}(a)^{-\frac{%
T_{3}n}{2}}
\end{eqnarray*}
and hence $d(x^{-n}u,V_{x^{-1}})\leq \Lambda _{k}(a)^{-kT_{3}}$ if $n\geq
n_{6}=n_{6}(l,d,(k_{i})_{i},(T_{i})_{i}).$

We check that similarly, $cx^{n}c^{-1}$ maps $B_{k}^{-}\cup B_{k}$into $%
B_{k}(c)$ and $cx^{-n}c^{-1}$ maps $B_{k}^{-}\cup B_{k}$ into $B_{k}^{-}(c)$
as soon as $n$ is larger that some fixed number depending only on the data $%
(l,d,(k_{i})_{i},(T_{i})_{i}).$

Finally we check that all balls $B_{k}^{-}$, $B_{k}$ $B_{k}^{-}(c)$, $%
B_{k}(c)$ are disjoint, since $d(V_{x},V_{x^{-1}})\geq d(V_{x},H_{x})\geq
\Lambda _{k}(a)^{-2T_{3}}$ and $d(V_{x},cV_{x^{-1}})\geq
d(H_{x},cV_{x^{-1}})\geq ||F||^{-T_{2}}\geq \Lambda _{k}(a)^{-k}$ and
similarly $d(cV_{x},V_{x^{-1}})\geq \Lambda _{k}(a)^{-k}$ and $%
d(cV_{x},cV_{x^{-1}})\geq ||F||^{-2dk_{5}}\Lambda _{k}(a)^{-2T_{3}}\geq
\Lambda _{k}(a)^{-(k-1)T_{3}}$.

It follows that $x^{n}$ and $cx^{n}c^{-1}$ play ping-pong on $\Bbb{P}%
(k^{d}), $ hence generate a free subgroup. This ends the proof of Lemma \ref
{pair}.

\section{Height bounds and proof of Theorem \ref{main}\label{finalproof}}

\subsection{A Product formula for subspaces\label{pfsubs}}

In this paragraph we define the \textit{adelic distance} $\delta (V;W)$
between two projective subspaces and we give a \textit{product formula} (\ref
{pfsub}) relating it to the Arakelov heights of $V,W$ and $V+W.$

In Paragraph \ref{FubStudSec} we recalled the Fubini-Study metric on $\Bbb{P}%
(k^{d})$, where $k$ is a local field. In particular, we had formula (\ref
{distsubspace}), which gives the distance between two projective linear
subspaces. If $K$ is a global field with prime field $K_{0}$ and $V$ and $W$
are disjoint projective linear subspace of $\Bbb{P}(K^{d})$, we can put
together the local distances (i.e. at each place of $K$) in a way similar to
the way the height of an algebraic number is defined. Namely we set: 
\begin{equation}
\delta (V;W)=\frac{1}{[K:K_{0}]}\sum_{v\in V_{K}}n_{v}\cdot \log \frac{1}{%
d_{v}(V,W)}  \label{hsub}
\end{equation}
where $d_{v}(\cdot ,\cdot )$ is the Fubini-Study metric on $\Bbb{P}%
(K_{v}^{d}).$ Each term in this sum is non negative. In fact, we see from (%
\ref{distsubspace}) that $\delta (V;W)$ is linked to the Arakelov heights
(see Paragraph \ref{Arak}) in the following simple way: 
\begin{equation}
\delta (V;W)=h_{Ar}(V)+h_{Ar}(W)-h_{Ar}(V+W)\leq h_{Ar}(V)+h_{Ar}(W)
\label{pfsub}
\end{equation}
This can be seen as a \textit{product formula for subspaces}, since when $V$
and $W$ are points in $\Bbb{P}^{1}(\overline{\Bbb{Q}})$ it reduces to the
classical product formula on $\overline{\Bbb{Q}}.$

Note moreover that we can similarly define $\delta (V^{\bot },W^{\bot })$
just as $\delta (V,W)$ in the projective space of the dual vector space $%
(K^{d})^{*}.$ Since $h_{Ar}(V)=h_{Ar}(V^{\bot })$ (see \cite{Bom}), we also
have 
\begin{equation*}
\delta (V^{\bot },W^{\bot })\leq h_{Ar}(V)+h_{Ar}(W)
\end{equation*}

We will often denote by $\delta _{v}(V;W)$ the term of the sum in (\ref{hsub}%
) corresponding to the place $v$, so that 
\begin{equation*}
\delta (V;W)=\frac{1}{[K:K_{0}]}\sum_{v\in V_{K}}n_{v}\cdot \delta _{v}(V;W)
\end{equation*}

\subsection{The Eskin-Mozes-Oh Escape Lemma\label{emo}}

In this paragraph we recall a crucial Lemma due Eskin-Mozes-Oh, which allows
to ``escape from algebraic subvarieties in bounded time''.

Recall Bezout's theorem about the intersection of finitely many algebraic
subvarieties (see for instance \cite{Sch}), namely:

\begin{thm}[Generalized Bezout theorem]
\label{thm:Bezout} Let $K$ be a field, and let $Y_{1},\ldots ,Y_{p}$ be pure
dimensional algebraic subvarieties of $K^{n}$. Denote by $W_{1},\ldots ,W_{q}
$ the irreducible components of $Y_{1}\cap \ldots \cap Y_{p}$. Then 
\begin{equation*}
\sum_{i=1}^{q}\text{deg}(W_{i})\leq \prod_{j=1}^{p}\text{deg}(Y_{j}).
\end{equation*}
\end{thm}

Let $K$ be a field and let $X$ be an algebraic variety over $K$. We set $%
s(X) $ to be the sum of the degree and the dimension of each of its
geometrically irreducible components. The following result was shown in \cite
{EMO}, Lemma 3.2:

\begin{lem}
\cite{EMO}\label{EMOLemma} Given an integer $m\geq 1$ there is $N=N(m)$ such
that for any field $K$, any integer $d\geq 1$, any $K$--algebraic subvariety 
$X$ in $GL_{d}(K)$ with $s(X)\leq m$ and any (not necessarily symmetric)
subset $F\subset GL_{d}(K)$ which contains the identity and generates a
subgroup which is not contained in $X(K)$, we have $F^{N}\nsubseteq X(K)$.
\end{lem}

\subsection{Irreducible representations of Chevalley groups\label{irredrep}}

In this paragraph we define the linear irreducible representations $(\rho
_{\alpha },E_{\alpha })$ which are the possible candidates for the
projective representation where we will play ping-pong. We also set a
particular basis in each $E_{\alpha }$, which we use to define the height $%
h(\rho _{\alpha }(g))$ and then show Lemma \ref{rephei}.

Let $\Bbb{G}$ be a Chevalley group of adjoint type and $\frak{g}$ its Lie
algebra with $\Bbb{Z}$-structure $\frak{g}_{\Bbb{Z}}$. Let $T$ be a maximal
torus and $\frak{t}$ the corresponding Cartan subalgebra in $\frak{g.}$ Let $%
\Lambda _{R}$ be the lattice of roots in the dual of $\frak{g}$ which we
identify with the space $X(T)$ of characters of $T.$ Let $\Lambda _{W}$ be
the lattice of weights. We fix a set of positive roots $\Phi ^{+}$ and
inside a base of simple roots $\Pi .$ Since $\Bbb{G}$ is of adjoint type, to
every dominant weight $\lambda \in \Lambda _{R},$ there correspond a finite
dimensional absolutely irreducible representation $E$ of $\Bbb{G}$. Let $%
\{\pi _{\alpha }\}_{\alpha \in \Pi }\subset \Lambda _{W}$ be the fundamental
weights. For each $\alpha \in \Pi ,$ there is a smallest integer $k_{\alpha
}\in \Bbb{N}$ such that $k_{\alpha }\pi _{\alpha }\in \Lambda _{R}.$ Let $%
\chi _{\alpha }=k_{\alpha }\pi _{\alpha }$ be the corresponding dominant
weight and $(\rho _{\alpha },E_{\alpha })$ the corresponding absolutely
irreducible representation of $\Bbb{G}$.

For background on Chevalley groups and their representations, see
Steinberg's notes \cite{Stein}. Let also $(\rho _{0},E_{0})$ be the adjoint
representation. According to \cite{Stein} Section 2 Theorem 2, given an
absolutely irreducible representation $(\rho ,E)$ of $\Bbb{G}$, one may find
in each $E$ a lattice $\Lambda $ invariant under the action of $\rho (\Bbb{G}%
(\Bbb{Z}))$ and a basis of $\Lambda $ which is made of weight vectors. Let
us choose this basis. It defines a standard norm $||\cdot ||_{k}$ on $%
\Lambda \otimes _{\Bbb{Z}}k$ and it also defines a height $h$ in $%
SL(E_{\alpha })$ as in Section \ref{normhei}. Recall that $\Omega $ is
either $\overline{\Bbb{Q}}$ or $\overline{\Bbb{F}(t)}$ and $\varepsilon
_{\Omega }=1$ in the first case $0$ otherwise. We have:

\begin{lemma}
\label{rephei}There exists a constant $C_{0}>0$ such that for every finite
subset $F\in \Bbb{G}(\Omega )$ and every $\alpha \in \Pi ,$ $h(\rho _{\alpha
}(F))\leq C_{0}\cdot (h(\rho _{0}(F))+\varepsilon _{\Omega }).$
\end{lemma}

Let $\chi _{\rho }$ be the heighest weight of $\rho $, which belongs to the
root lattice. Let $L$ be the maximal coefficient appearing in the
decomposition of $\chi _{\rho }$ as a sum of simple roots (let $L_{0}$ the
corresponding integer for $\rho _{0}=Ad$). Let $M$ be the smallest positive
integer such that $M\chi _{\rho }\geq \alpha $ for every $\alpha \in \Pi $
(for the order defined by $\Pi $). Then Lemma \ref{rephei} follows from:

\begin{lemma}
\label{normcomp}For every local field $k$, there is a constant $%
c_{0}=c_{0}(\rho ,k)>0$ such that for every $g\in \Bbb{G}(k),$ we have 
\begin{equation}
\frac{1}{c_{0}}||Ad(g)||_{k}^{\frac{1}{L_{0}M}}\leq ||\rho (g)||_{k}\leq
c_{0}||Ad(g)||_{k}^{L}  \label{normrepcomp}
\end{equation}
and $c_{0}=1$ unless $k$ is Archimedean.
\end{lemma}

\proof%
Let $\Bbb{K}_{k}=\Bbb{G}(\mathcal{O}_{k})$ when $k$ is ultrametric. Then $%
\rho (\Bbb{K}_{k})$ and $Ad(\Bbb{K}_{k})$ preserve the norm. By Cartan's $%
\Bbb{K}_{k}T\Bbb{K}_{k}$ decomposition, it suffices to prove the
inequalities for $g$ in the maximal torus $T.$ But then $||\rho
(g)||_{k}=|\chi _{\rho }(g)|_{k}$ and $\max \{|\alpha (g)|_{v},\alpha \in
\Pi \}^{\frac{1}{M}}\leq |\chi _{\rho }(g)|_{k}\leq \max \{|\alpha
(g)|_{v},\alpha \in \Pi \}^{L}.$ And $\max \{|\alpha (g)|_{v},\alpha \in \Pi
\}\leq ||Ad(g)||_{k}\leq \max \{|\alpha (g)|_{v},\alpha \in \Pi \}^{L_{0}}$.
Hence $(\ref{normrepcomp})$.

When $k$ is Archimedean, $\Bbb{K}_{k}$ stabilizes another norm $||\cdot
||_{k,new}$ on $\Lambda \otimes _{\Bbb{Z}}k$ (resp. $\Lambda _{R}\otimes _{%
\Bbb{Z}}k$). For this new norm the same argument gives $(\ref{normrepcomp}).$
Since the two norms are equivalent, this gives us the constant $c_{0}.$%
\endproof%

\subsection{Combined adelic distance}

In this paragraph, we define the combined adelic distance $\delta (F)=\delta
^{1}(F)+\delta ^{2}(F)$ of all adelic distances $\delta (V;W)$ where $V$ and 
$W$ range over the relevant projective subspaces involved in the ping-pong
conditions from Section \ref{PP}.

Let $K$ be a global field. Let $(q_{i})_{1\leq i\leq 5}$ be five positive
integers. Given $a\in \Bbb{G}(K)$ and $\alpha \in \Pi \cup \{0\},$ let $%
\mathcal{B}_{a,\alpha }$ be the set of elements $b\in \Bbb{G}(K)$ such that $%
^{t}\rho _{\alpha }(b)(V^{c})^{\bot }\nsubseteq W^{\bot }$ and $^{t}\rho
_{\alpha }(b^{-1})(V^{c})^{\bot }\nsubseteq W^{\bot }$ for every $(V,W)\in 
\mathcal{A}_{\alpha }(a),$ where $\mathcal{A}_{\alpha }(a)$ is the set of
couples $(V,W)$ of $\rho _{\alpha }(a)$-admissible (see def. \ref{admiss})
non-trivial linear subspaces of $E_{\alpha }$ such that $\dim (V)=1.$ Given $%
a,b\in \Bbb{G}(K)$ with $b\in \mathcal{B}_{a,\alpha },$ let $\mathcal{T}%
_{a,b,\alpha }$ be the set of elements $t\in \Bbb{G}(K)$ such that $\rho
_{\alpha }(t)V\nsubseteq W^{c}+W\cap b^{-1}V^{c}$ and $\rho _{\alpha
}(t^{-1})V\nsubseteq W^{c}+W\cap bV^{c}$ for every $V,W\in \mathcal{A}%
_{\alpha }(a)$ (note that since $b\in \mathcal{B}_{a,\alpha }$, $W^{c}+W\cap
b^{-1}V^{c}$ and $W^{c}+W\cap bV^{c}$ are hyperplanes).

Recall from Paragraph \ref{emo} that given an algebraic variety $\mathcal{Z}$
over the algebraically closed field $\Omega ,$ we denote by $s(\mathcal{Z})$
the sum of the dimension and degree of its irreducible components. Given two
non-trivial subspaces $V$ and $W$ in $E_{\alpha }$ the set of all $g\in
GL(E_{\alpha })$ such that $gW\subset V$ or $g^{-1}W\subset V$ is a Zariski
closed subset $\mathcal{Z}_{V,W}$ of $GL(E_{\alpha }).$ Moreover $s(\mathcal{%
Z}_{V,W})$ is bounded independently on $V$ and $W$ since the one can pass
from one $\mathcal{Z}_{V,W}$ to the other by multiplying on the left and
right by some automorphism in $GL(E_{\alpha }).$ From these remarks and
Lemma \ref{EMOLemma} we obtain:

\begin{lemma}
\label{nonempty}There is a positive integer $q_{0}$ such that for any field $%
K$ and any finite subset $F$ of $\Bbb{G}(K)$ containing $1$ and generating a
Zariski-dense subgroup, any $\alpha \in \Pi \cup \{0\}$ and any $a\in \Bbb{G}%
(K)$ and $b\in \mathcal{B}_{a,\alpha },$ the set $F^{q_{0}}$ intersects $%
\mathcal{B}_{a,\alpha }$ non trivially and the set $F^{q_{0}}$ intersects $%
\mathcal{T}_{a,b,\alpha }$ non trivially.
\end{lemma}

We now fix the values of $q_{2},q_{3}$ and $q_{5}$ to be equal to this $%
q_{0}.$ The values of $q_{1}$ and $q_{4}$ will be specified later. Let $%
\mathcal{Q}_{\alpha }$ be the set of $3$-tuples $(a,b,t)$ such that $a\in
F^{q_{1}},$ $b\in F^{q_{2}}\cap \mathcal{B}_{a,\alpha },$ and $t\in
F^{q_{3}}\cap \mathcal{T}_{a,b,\alpha }.$ Let $\mathcal{R}_{\alpha }$ be the
set of couples $(x,c)$ such that $x\in F^{q_{4}},$ $c\in F^{q_{5}}\cap 
\mathcal{B}_{x,\alpha }.$ Lemma \ref{nonempty} ensures that if $F$ generates
a Zariski-dense subgroup, then for any $a\in F^{q_{1}}$ there are $b,t$ such
that $(a,b,t)\in \mathcal{Q}_{\alpha }$ and also for any $x\in F^{q_{4}}$
there is $c$ such that $(x,c)\in \mathcal{R}_{\alpha }.$ Now define for any
finite symmetric subset $F$ in $\Bbb{G}(K),$ and $i=1,2$%
\begin{equation*}
\delta ^{i}(F)=\sum_{\alpha \in \Pi \cup \{0\}}\delta _{\alpha }^{i}(F)
\end{equation*}
and 
\begin{eqnarray*}
\delta _{\alpha }^{1}(F) &=&\sum_{(a,b,t)\in \mathcal{Q}_{\alpha }}\delta
_{\alpha ,(a,b,t)}(F) \\
\delta _{\alpha }^{2}(F) &=&\sum_{(x,c)\in \mathcal{R}_{\alpha }}\delta
_{\alpha ,(x,c)}(F)
\end{eqnarray*}
and 
\begin{eqnarray*}
\delta _{\alpha ,(a,b,t)}(F) &=&\sum_{(V,W)\in \mathcal{A}_{\alpha
}(a)}\delta (^{t}\rho _{\alpha }(b)(V^{c})^{\bot };W^{\bot })+\delta
(^{t}\rho _{\alpha }(b^{-1})(V^{c})^{\bot };W^{\bot })+ \\
&&\delta (\rho _{\alpha }(t)V;W^{c}+W\cap \rho _{\alpha
}(b^{-1})V^{c})+\delta (\rho _{\alpha }(t^{-1})V;W^{c}+W\cap \rho _{\alpha
}(b)V^{c})
\end{eqnarray*}
\begin{equation*}
\delta _{\alpha ,(x,c)}(F)=\sum_{(V,W)\in \mathcal{A}_{\alpha }(x)}\delta
(\rho _{\alpha }(c)V;W)+\delta (\rho _{\alpha }(c^{-1})V;W)
\end{equation*}

\subsection{Height bounds for subspace separation}

In this paragraph, applying the results of Paragraphs \ref{Arak} and \ref
{pfsubs}, we obtain $(\ref{delta1})$ and $(\ref{delta2})$ which give bounds
for the combined adelic distances $\delta ^{i}(F)$ in terms of the height $%
h(F)$ and the number of elements in $F$ only.

Namely, if $r=rank(\Bbb{G)}$ and $D=(r+1)\max_{\alpha \in \Pi \cup
\{0\}}24\cdot 4^{d_{\alpha }}d_{\alpha }^{2}$ with $d_{\alpha }=\dim
E_{\alpha },$ we have for every $(a,b,t)\in \mathcal{Q}_{\alpha }(F)$, 
\begin{eqnarray*}
\delta _{\alpha ,(a,b,t)}(F) &\leq &\sum_{V,W\in \mathcal{A}_{\alpha
}(a)}d_{\alpha }^{2}\cdot (h(\rho _{\alpha }(b))+h(\rho _{\alpha
}(t)))+4(h_{Ar}(W)+h_{Ar}(V))+ \\
&&+2(h_{Ar}(W^{c})+h_{Ar}(V^{c})) \\
&\leq &h(\rho _{\alpha }(F^{q_{0}}))+12\cdot 4^{d_{\alpha }}\cdot \max
\{h_{Ar}(W),W\text{~}a\text{-admissible}\} \\
&\leq &24\cdot 4^{d_{\alpha }}d_{\alpha }^{2}\cdot (q_{0}+q_{1})h(\rho
_{\alpha }(F))+12\cdot 4^{d_{\alpha }}d_{\alpha }^{2}\cdot \varepsilon
_{\Omega }\log 2
\end{eqnarray*}
Hence using Lemma \ref{rephei}, 
\begin{eqnarray*}
\delta ^{1}(F) &\leq &D|F|^{q_{1}+2q_{0}}(q_{0}+q_{1})\cdot h(\rho _{\alpha
}(F))+D\cdot \varepsilon _{\Omega }\log 2 \\
&\leq &D|F|^{q_{1}+2q_{0}}(q_{0}+q_{1})C_{0}\cdot (h(Ad(F))+\varepsilon
_{\Omega })
\end{eqnarray*}
Note that if $char(\Omega )=0,$ then by the Height Gap Theorem \ref{HeGa},
we have $h(Ad(F))\geq \widehat{h}(Ad(F))\geq g>0$ where $g$ is the gap. So
at any case for all characteristic, 
\begin{equation}
\delta ^{1}(F)\leq D_{1}\left( \frac{|F|}{5}\right) ^{q_{1}+2q_{0}}h(Ad(F)),
\label{delta1}
\end{equation}
where $D_{1}=D5^{q_{1}+2q_{0}}(q_{0}+q_{1})C_{0}(1+g^{-1}).$ Similarly one
obtains 
\begin{equation}
\delta ^{2}(F)\leq D_{2}\left( \frac{|F|}{5}\right) ^{q_{4}+q_{0}}h(Ad(F)),
\label{delta2}
\end{equation}
where $D_{2}=D5^{q_{4}+q_{0}}(q_{0}+q_{4})C_{0}(1+g^{-1}).$

\subsection{Proof of Theorem \ref{main}}

The proof is done in three steps. First we reduce to the situation when $F$
generates a Zariski-dense subgroup in $\Bbb{G}(\Omega )$ where $\Bbb{G}$ is
a simple Chevalley group of adjoint type to be chosen among a finite list of
such. Second we show that we may assume that $F=\{1,X,X^{-1},Y,Y^{-1}\}$,
i.e. $F$ is a symmetric set with $4$ elements plus the identity. And
finally, in the third and most difficult step, we check that there exists a
place $v$ of the field $K$ of coefficients for which the sufficient
conditions $(i)$ to $(vi)$ stated in Section \ref{PP} are fulfilled with
some explicit choice of constants depending only on $\Bbb{G}$, and thus
yield the desired ping-pong pair.

\begin{remark}
It is not clear whether or not the assumption $F$ symmetric is a necessary
condition in Theorem \ref{main}. Our proof however requires this assumption
(see Remark \ref{benoist}). If one needs only a free semi-group instead of a
free group, then it is not necessary.
\end{remark}

\subsubsection{Preliminary reductions. \newline
}

In this paragraph, we prove the first two steps, Claims 1 and 2. We have:

\textbf{Claim 1:} In Theorem \ref{main}, we may assume that $F$ generates a
Zariski-dense subgroup in $\Bbb{G}(\Omega )$ where $\Bbb{G}$ is a simple
Chevalley group of adjoint type.

\proof%
Since $F$ generates a non virtually solvable subgroup $\left\langle
F\right\rangle $, the connected component $\Bbb{G}^{0}$ of the
Zariski-closure $\Bbb{G}$ of $\left\langle F\right\rangle $ is not solvable.
Moding out by the solvable radical of $\Bbb{G}^{0},$ which is a normal
subgroup of $\Bbb{G}$, we see that we can assume that $\Bbb{G}^{0}$ is a non
trivial semisimple algebraic group. We let $\Bbb{G}$ act on $\Bbb{G}^{0}$ by
conjugation we obtain a homomorphism of $\Bbb{G}$ in $Aut(\Bbb{G}_{ad}^{0}%
\Bbb{)}$ where $\Bbb{G}_{ad}^{0}$ is the adjoint group of $\Bbb{G}^{0}$
whose image contains $\Bbb{G}_{ad}^{0}.$ However by \cite{Bor} IV.14.9. $Aut(%
\Bbb{G}_{ad}^{0}\Bbb{)}/\Bbb{G}_{ad}^{0}$ is a subgroup of the automorphisms
of the Dynkin diagram of $\Bbb{G}$. In particular it is a finite group whose
order is bounded in terms of $\dim \Bbb{G}$ only, hence in terms of $d$
only. Recall (see for instance \cite{HG} Lemma 4.6.),

\begin{lemma}
\label{finiteindex}Let $F$ be a finite subset of a group $\Gamma $
containing $1$. Assume that the elements of $F$ (together with their
inverses) generate $\Gamma .$ Let $\Gamma _{0}$ be a subgroup of index $k$
in $\Gamma .$ Then $F^{2k+1}$ contains a generating set of $\Gamma _{0}.$
\end{lemma}

Applying this lemma, we may therefore assume that $\Bbb{G}=\Bbb{G}_{ad}^{0}$
is a semisimple algebraic group of adjoint type. Further projecting to one
of the simple factors, we may assume that $\Bbb{G}$ is a simple algebraic
group of adjoint type over $\Omega .$ As $\Omega $ is algebraically closed, $%
\Bbb{G(}\Omega )$ is the group of $\Omega $-points of a Chevalley group (see 
\cite{Stein}). 
\endproof%

Let $\mathcal{O}$ be the Zariski-open subset of $\Bbb{G}\times \Bbb{G}$
obtained in Theorem \ref{GoodPos}.

\textbf{Claim 2:} In Theorem \ref{main}, we may assume that $%
F=\{1,X,X^{-1},Y,Y^{-1}\}$ for some $(X,Y)\in \mathcal{O}(\Omega ).$

\proof%
This claim was already proven in Proposition 4.14 of \cite{HG} in the
special case of characteristic $0$ making key use of Jordan's theorem about
finite subgroups of $GL_{n}(\Bbb{C}).$ This argument fails in positive
characteristic so we now give a different (and more involved) argument. Let $%
\Bbb{G}$ be a simple Chevalley group. Following an idea used in \cite{uti}
Section 7, we have:

\begin{lemma}
\label{zardense}Then there is a proper closed subvariety $\mathcal{W}$ of $%
\Bbb{G\times G}$ such that, for any choice of $\Omega $, every pair $%
(x,y)\notin \mathcal{W}(\Omega )$ with $x$ of infinite order generates a
Zariski-dense subgroup of $\Bbb{G}$.
\end{lemma}

\proof%
Let $\frak{g}$ be the Lie algebra of $\Bbb{G}$ (see \cite{Bor} I.3.5). Let $%
\mathcal{W}$ be the subset of pairs $(x,y)$ in $\Bbb{G}(\Omega )\Bbb{\times G%
}(\Omega )$ such that the associative subalgebra of $End(\frak{g})$
generated by $Ad(x)$ and $Ad(y)$ is proper. Note that $\mathcal{W}$ is a
closed algebraic subset with equations over $\Bbb{Z}$. It is also proper
because one can construct pairs $(x,y)$ for which the group they generate
acts irreducibly on $\frak{g}$ (see for instance \cite{Bourb} VIII. 2. ex.8.
and \cite{BaLa} \S 3). Suppose $(x,y)\notin \mathcal{W}(\Omega )$ and $x$
has infinite order. Let $\Bbb{H}$ be the Zariski closed subgroup generated
by $x$ and $y$. Then $\dim \Bbb{H}\geq 1$ and the Lie algebra of $\Bbb{H}$
is non trivial and invariant under $Ad(x)$ and $Ad(y)$, hence equal to $%
\frak{g}$. By \cite{Bor} I.3.6 we conclude that $\Bbb{H=G}$. 
\endproof%

In order to apply this lemma, we show:

\begin{lemma}
\label{pos}There is a constant $N=N(d)\in \Bbb{N}$ such that, for any choice
of $\Omega $, if $F$ is a finite symmetric subset of $\Bbb{G}(\Omega )$
containing $1$ and generating a Zariski dense subgroup, one may find a
subset $F_{0}$ of $F^{N}$ such that for all integers $n\geq 1$ the subset $%
F_{0}^{n}$ is made only of elements of infinite order and the subgroup
generated by $F_{0}$ and $F_{0}^{-1}$ is Zariski dense in $\Bbb{G}$.
\end{lemma}

Before going into the proof of Lemma \ref{pos} let us explain how we deduce
Claim 2 from this.

\textit{Proof of Claim 2}. By Lemma \ref{pos}, we can replace $F$ by $F_{0}$%
. Now according to Lemma \ref{EMOLemma} applied to $\Bbb{G\times G}$ and $%
F_{0}\times F_{0}$ there is a constant $M\in \Bbb{N}$ depending only on $%
\mathcal{W}$ and $\mathcal{O}$, hence on $d$ only, such that $F_{0}^{M}$
contains a pair $(x,y)$ such that $(x,y)\in \mathcal{O}$ and $(x,y)\notin 
\mathcal{W}$. By Lemma \ref{pos}, $x$ has infinite order, hence by Lemma \ref
{zardense}, $x$ and $y$ generate a Zariski dense subgroup of $\Bbb{G}$, and
Claim 2 is proved.

\proof[Proof of Lemma \ref{pos}]%
Let $d_{0}=\dim \Bbb{G}$. We have:

\begin{lemma}
\label{info}There is a constant $N_{0}=N_{0}(d_{0})\in \Bbb{N}$ and $k\leq
d_{0}$ elements $\alpha _{1},...,\alpha _{k}$ in $F^{N_{0}}$ of infinite
order and such that the connected components $C_{i}$ of the Zariski closures
of each cyclic subgroup generated by each $\alpha _{i}$ together generate $%
\Bbb{G}$ as an algebraic group.
\end{lemma}

\proof%
First we check that there is some $\alpha $ of infinite order in a bounded
power of $F,$ say $F^{N_{1}}$. This follows from Theorem \ref{HeGa} in
characteristic $0$ (see Corollary 1.2). In positive characteristic it
follows directly from the fact that $\left\langle F\right\rangle $ is finite
as soon as $\widehat{h}(F)=0$ (Lemma \ref{zerocase}) and Lemma \ref{NvsS} $%
(a)$ which says that $F^{d_{0}^{2}}$ already contains an element with
eigenvalue of absolute value $>1$.

Let $\left\langle \alpha \right\rangle $ be the cyclic group generated by $%
\alpha $ and $C_{1}$ the connected component of its Zariski closure. Then $%
\dim C_{1}=1.$ Set $\alpha _{1}=\alpha .$ Suppose $j\geq 1$ and we have
built $\alpha _{1},...,\alpha _{j}$ and let $C_{i}$ be the connected
component of the Zariski closure of $\left\langle \alpha _{i}\right\rangle $
and $\Bbb{H}_{i}$ the algebraic subgroup generated by all $C_{m}$ for $1\leq
m\leq i.$ We show by induction that $\alpha _{i}=w_{i-1}\alpha w_{i-1}^{-1}$
for some $w_{i-1}\in F^{i-1}$ and $\dim \Bbb{H}_{i}\geq i.$ If $\Bbb{H}%
_{j}\neq \Bbb{G}$, as $\Bbb{G}$ is simple and $\left\langle F\right\rangle $
Zariski dense, there must exist some $\beta _{j}\in F$ such that $\beta _{j}%
\Bbb{H}_{j}\beta _{j}^{-1}\neq \Bbb{H}_{j}$, hence some $i\leq j$ such that $%
\beta _{j}C_{i}\beta _{j}^{-1}$ is not contained in $\Bbb{H}_{j}.$ Let $%
\alpha _{j+1}=\beta _{j}\Bbb{\alpha }_{i}\beta _{j}^{-1}$, i.e. $w_{j}=\beta
_{j}w_{i-1}\in F^{j}.$ Then $C_{j+1}=\beta _{j}C_{i}\beta _{j}^{-1}$ and $%
\dim \Bbb{H}_{j+1}\geq \dim \Bbb{H}_{j}+1.$ 
\endproof%

We look at $\Bbb{G}$ viewed inside $SL(\frak{g})$ via the adjoint
representation. We know from Theorem \ref{HeGa} and Lemma \ref{NvsS} that
either there is a non archimedean place $v$ of $\Omega $ for which $\Lambda
_{v}(F^{d_{0}^{2}})>1$ or there is an archimedean place for which $\Lambda
_{v}(F^{d_{0}^{2}})>1+\varepsilon $ where $\varepsilon $ is the Height Gap.
Let $f\in F^{d_{0}^{2}}$ be such that $\Lambda _{v}(f)=\Lambda
_{v}(F^{d_{0}^{2}}).$ At any case, in one of boundedly many irreducible
representations of $\Bbb{G}$ over the local field $K_{v},$ $f$ acts as a
proximal transformation with a contracting eigenvalue and its action on the
associated projective space $\Bbb{P}(K_{v}^{D})$ is described by Lemma \ref
{prox}. Let $v_{f}$ be its attracting point and $H_{f}$ be the repelling
hyperplane. By Lemma \ref{NvsS}, we may conjugate $F$ inside $GL_{D}(K_{v})$
so that $||F||_{v}$ is less than say $\Lambda _{v}(f)^{c_{0}}$ where $c_{0}$
is some constant depending only on $d_{0}.$ Up to changing $f$ into $%
f^{c_{0}}$ we may assume that $||F||_{v}\leq \Lambda _{v}(f).$

Let $\alpha _{1},...,\alpha _{k}$ be the elements from Lemma \ref{info}.
According to Lemma \ref{Cayley}, there is some $n_{i}\in [1,d_{0}]$ such
that $d(\alpha _{i}^{n_{i}}v_{f},H_{f})^{-1}$ is bounded above by some
bounded power of $||F||_{v}.$ If follows from Lemma \ref{prox} that there is
compact subset $C$ of the projective space $\Bbb{P}(K_{v}^{D})$ which is the
complement of some neighborhood of $H_{f},$ such that after replacing $f$ by
some bounded power of it if necessary, the elements $f,$ $\alpha
_{1}^{n_{1}}f,$ ..., $\alpha _{k}^{n_{k}}f$ are all proximal, send $C$
inside itself, and have a Lipschitz constant $<1$ on $C$. Let $%
F_{0}=\{f,\alpha _{1}^{n_{1}}f,...,\alpha _{k}^{n_{k}}f\}.$ We check that $%
F_{0}$ satisfies the desired conditions. It lies in a bounded power of $F$,
every positive word with letters in $F_{0}$ preserves $C$ and is proximal by
Tits' converse Lemma \ref{TitsConv}, hence of infinite order. Finally the
group $\left\langle F_{0}\right\rangle $ generated by $F_{0}$ contains each $%
\left\langle \alpha _{i}^{n_{i}}\right\rangle $, hence its Zariski closure
contains the connected component $C_{i}$ of the cyclic group $\left\langle
\alpha _{i}\right\rangle .$ Since the $C_{i}$'s generate $\Bbb{G}$ as an
algebraic group by Lemma \ref{info}, we get that $\left\langle
F_{0}\right\rangle $ is Zariski dense, and this ends the proof of Lemma \ref
{pos}. 
\endproof%

\subsubsection{End of the proof of Theorem \ref{main}. \newline
}

So from now on we assume that $\Bbb{G}$ is a simple Chevalley group of
adjoint type over $\Omega $ viewed as embedded inside $SL(\frak{g})$ ($\frak{%
g}=Lie(\Bbb{G)}$) where it acts via the adjoint representation. We also
assume that $F=\{1,X,X^{-1},Y,Y^{-1}\}$ generates a Zariski-dense subgroup
of $\Bbb{G}$ and $(X,Y)$ lies in the Zariski-open subset $\mathcal{O}$
defined in Theorem \ref{GoodPos}.

\paragraph{\textbf{Constants. }\newline
}

We now define or recall our constants. All these constants depend only on $%
\Bbb{G}$ (equivalently only on $\dim \Bbb{G}$) and not on the field of
coefficients we choose. And this is all that matters, so the reader may
freely ignore their precise definition, all the more so since we did not try
at all to give the best constants we could. However there dependence and
order in which they are defined are important in the logic of the proof.

Recall that the constant $C_{k,1}$ from Section \ref{PP} was defined to be $%
1 $ if $k$ is ultrametric and equal to the dimension of the vector space if $%
k$ is Archimedean.

Below we set the value of $d$ to be the $\max d_{\alpha }$ where $d_{\alpha
}=\dim E_{\alpha }$ for $\alpha \in \Pi \cup \{0\}$ (recall that we chose to
denote by $E_{0}$ the adjoint representation, so $d_{0}=\dim \Bbb{G}$).

$D=(rk(\Bbb{G)}+1)\max_{\alpha \in \Pi \cup \{0\}}24\cdot 4^{d_{\alpha
}}d_{\alpha }^{2}.$

$L$ is defined to be the maximum coefficient in the expression of the
heighest weight $\chi _{\alpha }$ (for each $\alpha \in \Pi \cup \{0\}$) as
a sum of simple roots.

$M$ is the smallest positive integer $k$ such that $k\chi _{\rho _{\alpha
}}-\beta $ is positive for any choice of simple roots $\alpha ,\beta \in \Pi
.$

$c_{0}(v)$ is the maximum of the constants denoted $c_{0}$ in Lemma \ref
{normcomp} for each $\rho _{\alpha },$ $\alpha \in \Pi \cup \{0\}$ for a
given place $v$ ($c_{0}=1$ when $v$ is finite and $c_{0}(v)$ is a fixed
constant $c_{0}(\infty )$ if $v$ is infinite).

$c(d_{0})_{v}$ is the constant appearing in the Comparison Lemma, Lemma \ref
{NvsS}. It is $1$ if $v$ is finite, a fixed constant $c(d_{0})_{\infty }$ if 
$v$ is infinite.

$g$ is the Height Gap from Theorem \ref{HeGa} in $SL_{d_{0}}(\overline{\Bbb{Q%
}}).$

If $char(\Omega )>1,$ then we set $n_{1}=1,$ otherwise we set $n_{1}$ to be
the first integer such that $\exp (\frac{g}{4}\sqrt{\frac{n_{1}}{8d}})\geq
\max \{c(d_{0})^{-2},d^{4LMd},c_{0}(\infty )^{2LM}\}.$

$q_{0}$ is the integer obtained by escape in Lemma \ref{nonempty}.

$q_{1}$ is $n_{1}d^{2}.$

$\varepsilon _{0}$ is $\frac{1}{L}.$

$\varepsilon $ is $\varepsilon _{0}/12d^{2}.$

$T_{1}$ is the maximum of the integers $\tau _{1}(d_{\alpha },\varepsilon )$
obtained in Lemma \ref{veryprox} for each representation $\rho _{\alpha },$ $%
\alpha \in \Pi \cup \{0\}.$

$C$ is the constant from Theorem \ref{GoodPos}, applied to $\Bbb{G}$ inside $%
SL_{d_{0}}.$

$C_{0}$ is the constant from Lemma \ref{rephei}.

Let $m=48T_{1}n_{1}CL^{2}.$

Let $D_{1}=D5^{q_{1}+2q_{0}}(q_{0}+q_{1})C_{0}(1+g^{-1})$

Let $T_{0}=24CD_{1}LM$

Let $k_{1}=d^{2}m,$ $k_{2}=k_{3}=k_{5}=q_{0}.$

Let $T_{3}$ be the maximum of the integers $\tau _{3}$ obtained in Lemma \ref
{veryprox} for the above values of $d_{\alpha
},k_{1},k_{2},k_{3},\varepsilon _{0},\varepsilon ,T_{0}$ and $T_{1}$ for
each representation $\rho _{\alpha },$ $\alpha \in \Pi \cup \{0\}.$

Let $l$ be the maximum of the integers $l$ obtained in Lemma \ref{veryprox}
for the above values of $d_{\alpha },k_{1},k_{2},k_{3},\varepsilon
_{0},\varepsilon ,T_{0},T_{1}$ and $T_{3}$ for each representation $\rho
_{\alpha },$ $\alpha \in \Pi \cup \{0\}.$

Let $k_{4}=2k_{1}l+k_{2}+k_{3}.$

Let $q_{4}=n_{1}k_{4}.$

Let $D_{2}=D5^{q_{4}+q_{0}}(q_{4}+q_{0})C_{0}(1+g^{-1}).$

Let $T_{2}=24CD_{2}LM.$

Let $l_{2}$ be the maximum of each value $l_{2}(d_{\alpha },(k_{i})_{1\leq
i\leq 5},\varepsilon ,\varepsilon _{0},(T_{i})_{0\leq i\leq 3})$ obtained in
Lemma \ref{pair} for each representation $\rho _{\alpha },$ $\alpha \in \Pi
\cup \{0\}.$

\paragraph{\textbf{Choice of a place }$v$. \newline
}

Applying Theorem \ref{GoodPos}, we may change $F$ into a conjugate of it by
some element in $\Bbb{G}(\Omega )$ and hence get, summing $(\ref{RPos})$, $(%
\ref{delta1})$ and $(\ref{delta2}),$%
\begin{equation}
h(Ad(F))+\frac{1}{D_{1}}\delta ^{1}(F)+\frac{1}{D_{2}}\delta ^{2}(F)\leq
3C\cdot e(Ad(F))  \label{hb}
\end{equation}
Let $K$ be the (global) field generated by the coefficients of $F$.

\textbf{Claim}: There is a place $v$ of $K$ such that the following holds:\ $%
e_{v}>0$ if $v$ is a finite place, $e_{v}\geq \frac{g}{4}$ if $v$ is
infinite and in both cases

\begin{eqnarray}
\log ||Ad(F)||_{v} &\leq &12C\cdot e_{v}  \label{place} \\
\delta ^{1}(F)_{v} &\leq &12CD_{1}\cdot e_{v}  \notag \\
\delta ^{2}(F)_{v} &\leq &12CD_{2}\cdot e_{v},  \notag
\end{eqnarray}
where $e_{v}=\log E_{v}(Ad(F))$ and $\delta ^{i}(F)_{v}$ is the part of $%
\delta ^{i}(F)$ associated to $v,$ i.e. 
\begin{equation*}
\delta ^{i}(F)=\frac{1}{[K:K_{0}\Bbb{]}}\sum_{v\in V_{K}}n_{v}\cdot \delta
^{i}(F)_{v}
\end{equation*}

\textit{Proof of claim:} This is an easy verification. Indeed, splitting the
infinite part and the finite part write $e=e(Ad(F))=e_{\infty }+e_{f}$. If $%
e_{\infty }<\frac{e}{2},$ then $e\leq 2e_{f}$ and $(\ref{hb})$ implies 
\begin{equation*}
h_{f}(Ad(F))+\frac{1}{D_{1}}\delta _{f}^{1}(F)+\frac{1}{D_{2}}\delta
_{f}^{2}(F)\leq 6C\cdot e_{f}(Ad(F))
\end{equation*}
where the subscript $f$ means that we have restricted the sum to the finite
places. Then the existence of a finite place $v$ such that $e_{v}>0$ and $(%
\ref{place})$ holds is guaranteed. On the other hand, if $e_{\infty }\geq 
\frac{e}{2},$ then we have 
\begin{equation*}
h_{\infty }(Ad(F))+\frac{1}{D_{1}}\delta _{\infty }^{1}(F)+\frac{1}{D_{2}}%
\delta _{\infty }^{2}(F)\leq 6C\cdot e_{\infty }(Ad(F))
\end{equation*}
Let $V^{+}$ be the set of places $v\in V_{\infty }$ for which $e_{v}\geq 
\frac{e_{\infty }}{2}.$ We have 
\begin{equation*}
\frac{e_{\infty }}{2}\leq \frac{1}{[K:K_{0}]}\sum_{v\in V^{+}}n_{v}e_{v}
\end{equation*}
And 
\begin{equation*}
h_{\infty }(Ad(F))+\frac{1}{D_{1}}\delta _{\infty }^{1}(F)+\frac{1}{D_{2}}%
\delta _{\infty }^{2}(F)\leq 12C\cdot \frac{1}{[K:K_{0}]}\sum_{v\in
V^{+}}n_{v}e_{v}
\end{equation*}
which surely guarantees the existence of a place $v\in V^{+}$ such that $(%
\ref{place})$ holds, and as $v\in V^{+},$ $e_{v}\geq \frac{e}{4}\geq \frac{g%
}{4}.$ qed.

\paragraph{\textbf{Verification of the Ping-Pong conditions }$(i)$\textbf{\
to }$(vi)$\textbf{\ from Section \ref{PP}. }\newline
}

We are going to build an element $a\in F^{q_{1}}$ and choose an $\alpha \in
\Pi $ for which all the six conditions of Section \ref{PP} are going to be
satisfied with $a^{m}$ in place of $a$ and $\rho _{\alpha }(F^{n_{1}})$ in
place of $F.$

According to Lemma \ref{gr}, 
\begin{eqnarray*}
E_{v}(Ad(F^{n_{1}})) &\geq &E_{v}(Ad(F))^{\sqrt{\frac{n_{1}}{8d}}} \\
&\geq &\max \{c(d_{0})_{v}^{-2},C_{K_{v},1}^{4LMd},c_{0}(v)^{2LM}\}.
\end{eqnarray*}
From Lemma \ref{NvsS}, 
\begin{eqnarray*}
\Lambda _{v}(AdF^{q_{1}}) &\geq &\Lambda _{v}(AdF^{n_{1}d_{0}^{2}})\geq
c(d_{0})E_{v}(AdF^{n_{1}}) \\
&\geq &E_{v}(AdF^{n_{1}})^{\frac{1}{2}}\geq ||AdF||_{v}^{\frac{1}{12C}}>1
\end{eqnarray*}
Let us choose $a\in F^{q_{1}}$ such that $\Lambda _{v}(Ad(a))=\Lambda
_{v}(AdF^{q_{1}}).$ Let also $\alpha \in \Pi $ be such that $|\alpha
(a)|_{v}=\max \{|\beta (a)|_{v},\beta \in \Pi \}.$ Then the representation $%
\rho _{\alpha }$ (defined in Subsection \ref{irredrep}) and $Ad$ satisfy $%
\Lambda _{v}(Ad(a))\leq |\alpha (a)|_{v}^{L}$ and $\Lambda _{v}(\rho
_{\alpha }(a))\leq |\alpha (a)|_{v}^{L}$ by definition of $L.$ Hence $%
|\alpha (a)|_{v}>1.$ Moreover, by definition of $\rho _{\alpha }$ and $%
\varepsilon _{0}=\frac{1}{L}$ we have 
\begin{equation}
\frac{\Lambda _{v}(\rho _{\alpha }(a))}{\lambda _{v}(\rho _{\alpha }(a))}%
=|\alpha (a)|_{v}>1  \label{i}
\end{equation}
It follows that $\rho _{\alpha }(a)$ is proximal. Moreover 
\begin{equation}
\left( \frac{\Lambda _{v}(\rho _{\alpha }(a))}{\lambda _{v}(\rho _{\alpha
}(a))}\right) ^{\frac{1}{\varepsilon _{0}}}=|\alpha (a)|_{v}^{\frac{1}{%
\varepsilon _{0}}}\geq \Lambda _{v}(\rho _{\alpha }(a))  \label{iiia}
\end{equation}
and 
\begin{equation*}
\Lambda _{v}(\rho _{\alpha }(a))\geq |\alpha (a)|_{v}\geq \Lambda
_{v}(Ad(a))^{\frac{1}{L}}\geq ||AdF||_{v}^{\frac{1}{12CL}}
\end{equation*}
On the other hand, by Lemma \ref{normcomp}, we have $||\rho _{\alpha
}(F^{n_{1}})||_{v}\leq c_{0}(v)||AdF^{n_{1}}||_{v}^{L}\leq
||AdF^{n_{1}}||_{v}^{L+1}\leq ||AdF||_{v}^{n_{1}2L}$ so 
\begin{equation*}
\Lambda _{v}(\rho _{\alpha }(a))\geq ||\rho _{\alpha }(F^{n_{1}})||_{v}^{%
\frac{1}{24CL^{2}n_{1}}}
\end{equation*}
And 
\begin{equation}
\Lambda _{v}(\rho _{\alpha }(a^{m}))\geq ||\rho _{\alpha
}(F^{n_{1}})||_{v}^{T_{1}}  \label{iiib}
\end{equation}
Raising $a$ to the power $m,$ $(\ref{i})$ gives condition $(i),$ while $(\ref
{iiia})$ and $(\ref{iiib})$ give condition $(iii).$ On the other hand Lemma 
\ref{normcomp} gives 
\begin{equation}
||\rho _{\alpha }(F^{n_{1}})||_{v}\geq c_{0}(v)^{-1}||AdF^{n_{1}}||_{v}^{%
\frac{1}{LM}}\geq ||AdF^{n_{1}}||_{v}^{\frac{1}{2LM}}\geq C_{K_{v},1}^{2d}.
\label{lll}
\end{equation}
Hence condition $(ii)$ is fulfilled.

We now check $(iv)$ and $(v).$ By $(\ref{place})$ we have $\delta
^{1}(F)_{v}\leq 12CD_{1}\cdot e_{v}.$ Since $\delta ^{1}(F)_{v}$ is a sum of
positive terms, we get in particular for any $(b,t)$ such that $(a,b,t)\in 
\mathcal{Q}_{\alpha }$ (just pick one!) 
\begin{equation*}
\sum_{W}\delta (^{t}\rho _{\alpha }(b)H_{a}^{\bot };W^{\bot })_{v}+\delta
(^{t}\rho _{\alpha }(b^{-1})H_{a}^{\bot };W^{\bot })_{v}\leq 12CD_{1}\cdot
e_{v}
\end{equation*}
where $H_{a}$ the generalized eigenspace of $\rho _{\alpha }(a)$
corresponding to eigenvalues that are $<\Lambda _{v}(\rho _{\alpha }(a))$
(it is a hyperplane since $\rho _{\alpha }(a)$ is proximal), and the sum is
made over all non trivial $\rho _{\alpha }(a)$-admissible subspaces $W$.
This gives 
\begin{eqnarray*}
d_{v}(^{t}\rho _{\alpha }(b)H_{a}^{\bot };W^{\bot }) &\geq
&E_{v}(Ad(F^{n_{1}}))^{-12CD_{1}}\geq ||Ad(F^{n_{1}})||_{v}^{-12CD_{1}} \\
&\geq &||\rho _{\alpha }(F^{n_{1}})||_{v}^{-24CD_{1}LM} \\
&\geq &||\rho _{\alpha }(F^{n_{1}})||_{v}^{-T_{0}}
\end{eqnarray*}
Similarly 
\begin{equation*}
d_{v}(^{t}\rho _{\alpha }(b^{-1})H_{a}^{\bot };W^{\bot })\geq ||\rho
_{\alpha }(F^{n_{1}})||_{v}^{-T_{0}}
\end{equation*}
This proves $(iv).$ Condition $(v)$ is derived in exactly the same way.

Therefore we are in the situation where we may apply Lemma \ref{veryprox}.
It yields an element $x\in F^{n_{1}k_{4}}=F^{q_{4}}$ such that $\rho
_{\alpha }(x)$ is very proximal and satisfies the conclusions of Lemma \ref
{veryprox}. Pick $c\in F^{q_{0}}$ such that $(x,c)\in \mathcal{R}_{\alpha }$
(there are such $c$ by Lemma \ref{nonempty}). The third inequality in $(\ref
{place})$ gives for every $\rho _{\alpha }(x)$-admissible subspaces $V$ and $%
W$ with $\dim V=1,$ 
\begin{equation*}
\delta (\rho _{\alpha }(c)V;W)_{v}+\delta (\rho _{\alpha
}(c^{-1})V;W)_{v}\leq 12CD_{2}\cdot e_{v}
\end{equation*}
We may take $V=V_{x}$ or $V_{x^{-1}}$ and $W=H_{x}$ or $H_{x^{-1}}$ and this
indeed gives condition $(v)$ with $T_{2}=24CD_{2}LM.$

Finally Lemma \ref{pair} yields that $x^{n}$ and $cx^{n}c^{-1}$ generate a
free subgroup as soon as $n$ is larger than the constant $l_{2}$ ($l_{2}$ is
expressible explicitly in terms of all the other constants introduced so
far).

This ends the proof of the main theorem. Q.E.D.

\section{Applications\label{Corrrr}}

In this section we briefly discuss the corollaries. We shall be brief as
each of them is derived in exactly the same way as in the $GL_{2}$ case, so
we will refer the reader to the paper \cite{BGL2} for details. The proofs of
Corollaries \ref{poschar}, \ref{wdioph} and \ref{weq} rely only on the
characteristic $0$ part of Theorem \ref{main} and on a reformulation of that
theorem in terms of algebraic varieties. So we will content ourselves to
give this reformulation and briefly explain below what makes this
translation possible. The following fact is standard,

\begin{proposition}
\label{algsol}(see e.g. \cite{HG} Proposition 7.4.) Let $G=GL_{d}(\Bbb{C}).$
For every integer $k,$ let $\mathcal{V}$ be the set of $k$-tuples $%
(a_{1},...,a_{k})\in G^{k}$ which generate a virtually solvable subgroup.
Then $\mathcal{V}$ is a closed algebraic subvariety of $G^{k}.$
\end{proposition}

It is proved via the following proposition:

\begin{proposition}
There exists $N=N(d)$ such that $(a_{1},...,a_{k})\in G^{k}$ generates a
virtually solvable subgroup if and only if they leave invariant a common
finite subset of at most $N$ points on the flag variety $G/B$, where $B$ is
the subgroup of upper triangular matrices.
\end{proposition}

Let $N=N(d)$ be the integer obtained in the statement of Theorem \ref{main}
and let $B(n)$ be the ball of radius $n$ in the free group $F_{2}$ on two
generators. For $n\geq 1$ let $\mathcal{W}_{n}$ be the set of couples $%
(A,B)\in GL_{d}(\Bbb{C})^{2}$ such that for any words $w_{1}$ and $w_{2}$ in 
$B(N)$ there exists a word $w\in B(n)\backslash \{1\}$ such that $%
w(w_{1}(A,B),w_{2}(A,B))=1.$ Clearly $\mathcal{W}_{n}$ is a closed
subvariety of $GL_{d}(\Bbb{C})^{2}.$ We obtain:

\begin{proposition}
\label{equiv}Theorem \ref{main} for $K=\Bbb{C}$ is equivalent to the
statement: $\mathcal{W}_{n}\subset \mathcal{V}$ for every $n\geq 1$.
\end{proposition}

This allows to use the following effective version of Hilbert's
Nullstellensatz:

\begin{theorem}
\label{MasWus}(\cite{MW}) Let $r,d\in \Bbb{N}$, $h>0$ and $f,q_{1},...,q_{k}$
be polynomials in $\Bbb{Z}[X_{1},...,X_{r}]$ with logarithmic height at most 
$h$ and degree at most $d$. Assume that $f$ vanishes at all common zeros (if
any) of $q_{1},...,q_{k}$ in $\Bbb{C}[X_{1},...,X_{r}].$ Then there exist $%
a,e\in \Bbb{N}$ and polynomials $b_{1},...,b_{k}\in \Bbb{Z}[X_{1},...,X_{r}]$
such that 
\begin{equation*}
af^{e}=b_{1}q_{1}+...+b_{k}q_{k}
\end{equation*}
with $e\leq (8d)^{2^{r}}$, the total degree of each $b_{i}$ at most $%
(8d)^{2^{r}+1}$ and the logarithmic height of each $b_{i}$ as well as $a$ is
at most $(8d)^{2^{r+1}+1}(h+8d\log (8d)).$
\end{theorem}

Since the polynomial equations defining $\mathcal{W}_{n}$ have degree linear
in $n$ and height exponential in $n$, one can get from Theorem \ref{MasWus}
the desired bound on the degree and height of the $b_{i}$'s and on $a$ and $e
$. This readily allows to deduce Corollaries \ref{poschar}, \ref{wdioph} and 
\ref{weq} from Theorem \ref{main} and Corollary \ref{cogrowth}. Corollary 
\ref{exp} is derived in a similar fashion. See \cite{BGL2} for more details.

\begin{Ack}
I am grateful to J. Tits for his encouraging remarks at an early stage of
this project. I also thank J-F. Quint for telling me about his results and
those of Y. Benoist on proximal maps over the $p$-adics.
\end{Ack}

\end{document}